\documentclass[11pt]{article}

\usepackage{amssymb,latexsym}
\usepackage{epsfig}
\usepackage{eufrak}
\usepackage{amsmath}
\usepackage{mathrsfs}
\usepackage{color}

\newtheorem{theorem}{Theorem}
\newtheorem{lemma}{Lemma}

\newcommand{\be}{\begin{equation}}
\newcommand{\ee}{\end{equation}}
\newcommand{\bee}{\begin{eqnarray*}}
\newcommand{\eee}{\end{eqnarray*}}
\newcommand{\bel}{\begin{eqnarray}}
\newcommand{\eel}{\end{eqnarray}}
\newcommand{\bec}{\begin{cases}}
\newcommand{\eec}{\end{cases}}
\newcommand{\bem}{\begin{bmatrix}}
\newcommand{\eem}{\end{bmatrix}}

\newcommand{\la}{\label}
\newcommand{\li}{\left}
\newcommand{\ri}{\right}

\newcommand{\ovl}{\overline}
\newcommand{\udl}{\underline}

\newcommand{\lc}{\lceil}
\newcommand{\rc}{\rceil}
\newcommand{\lf}{\lfloor}
\newcommand{\rf}{\rfloor}

\newcommand{\ep}{\epsilon}
\newcommand{\vep}{\varepsilon}
\newcommand{\lm}{\lambda}

\newcommand{\si}{\sigma}

\newcommand{\de}{\delta}

\newcommand{\vpi}{\varpi}
\newcommand{\ga}{\gamma}
\newcommand{\Ga}{\Gamma}
\newcommand{\vse}{\vartheta}
\newcommand{\se}{\theta}

\newcommand{\ze}{\zeta}

\newcommand{\om}{\omega}
\newcommand{\Om}{\Omega}

\newcommand{\f}{\frac}

\newcommand{\cd}{\cdots}

\newcommand{\qu}{\quad}
\newcommand{\qqu}{\qquad}
\newcommand{\fa}{\forall}

\newcommand{\mscr}{\mathscr}
\newcommand{\mcal}{\mathcal}
\newcommand{\mbf}{\mathbf}
\newcommand{\bb}{\mathbb}

\newcommand{\wh}{\widehat}

\newcommand{\mrm}{\mathrm}
\newcommand{\bs}{\boldsymbol}

\newcommand{\LRA}{\Longleftrightarrow}
\newcommand{\sh}{\slash}

\newcommand{\tx}{\text}

\newcommand{\iy}{\infty}

\newcommand{\pa}{\partial}

\newcommand{\bed}{\begin{description}}
\newcommand{\eed}{\end{description}}
\newcommand{\bei}{\begin{itemize}}
\newcommand{\eei}{\end{itemize}}
\newcommand{\ben}{\begin{enumerate}}
\newcommand{\een}{\end{enumerate}}
\newcommand{\bib}{\bibitem}
\newcommand{\beL}{\begin{lemma}}
\newcommand{\eeL}{\end{lemma}}
\newcommand{\beT}{\begin{theorem}}
\newcommand{\eeT}{\end{theorem}}
\newcommand{\sect}{\section}

\newcommand{\bpf}{\begin{pf}}
\newcommand{\epf}{\end{pf}}
\newcommand{\bsk}{\bigskip}
\newcommand{\bi}{\binom}

\newcommand{\pfbox}{\hfill\mbox{$\Box$}}

\newenvironment{pf}{\paragraph*{Proof{\rm.}}}{\pfbox\bigskip}

\begin{document}

\title{{\bf A Theory of Truncated Inverse Sampling}
\thanks{The author had been previously working with Louisiana State University at Baton Rouge, LA 70803, USA,
and is now with Department of Electrical Engineering, Southern
University and A\&M College, Baton Rouge, LA 70813, USA; Email:
chenxinjia@gmail.com}}

\author{Xinjia Chen}

\date{November 2008}

\maketitle

\begin{abstract}

In this paper, we have established a new framework of truncated
inverse sampling for estimating mean values of non-negative random
variables such as binomial, Poisson, hyper-geometrical, and bounded
variables. We have derived explicit formulas and computational
methods for designing sampling schemes to ensure prescribed levels
of  precision and confidence for point estimators. Moreover, we have
developed interval estimation methods.

\end{abstract}

\sect{Introduction}

Parametric estimation based on sampling is an important branch of
mathematical statistics with ubiquitous applications across many
fields, from operation research, biology and medical science,
agriculture science, computer science, social science,
telecommunication engineering, control engineering, to name a few. A
wide class estimation problems of both theoretical and practical
significance can be put into the setting of estimating the mean
value of a random variable via sampling.  Familiar examples include
the estimation of binomial parameters, Poisson parameters, finite
population proportion, the mean of a bounded variable, and so on.  A
simple yet frequently used sampling scheme for estimating the mean
value a random variable $X$ is to draw samples of $X$ until the
sample sum is no less than a prescribed threshold and then take the
empirical mean as an estimate for the true mean value. This sampling
scheme, referred to as {\it inverse sampling},  was first studied by
Haldane \cite{H, H2} in the context of estimating a binomial
parameter.  Recently, inverse sampling  has been studied by Chen
\cite{Chen1, Chen2}, Dagum et al. \cite{Dagum} and Cheng
\cite{Cheng} for estimation of the mean of a bounded variable. Mendo
and Hernando \cite{MH} have revisited inverse sampling for
estimating binomial parameters.

Theoretically, there is no limit on the number of samples for
inverse sampling.  However, the practical situation is quite
contrary.  Due to the limitation of resources, almost every
practitioner would specify a maximum sample size on the sampling.
This means that the frequently used method is actually the {\it
truncated inverse sampling} scheme in the sense that sampling is
continued until the sample sum is no less than a prescribed
threshold or the number of samples reach the maximum sample size.

While the ideal inverse sampling has drawn extensive research
effort, little attention  has been paid to the theoretical issues of
the truly useful truncated inverse sampling scheme.  In this paper,
we shall investigate the essential theory of truncated inverse
sampling with a prevailing theme of error control.  We have answered
two equally central problems regarding pre-experimental planning and
post-experimental analysis. The first problem is on the
determination of the threshold value and the maximum sample size for
guaranteeing prescribed levels of precision and confidence of an
estimator.  The second problem is on interval estimation of the
parameter based on the observed data when the truncated inverse
sampling is completed.

The remainder of the paper is organized as follows. In Section 2, we
present our general results for truncated inverse sampling. In
Section 3, we consider the problem of estimating binomial
parameters. In Section 4, we discuss the estimation of the
proportion of a finite proportion. In Section 5, we discuss the
estimation of Poisson parameters. The estimation of the mean of a
bounded variable is investigated in Section 6.  Section 7 is the
conclusion.

Throughout this paper, we shall use the following notations. The
expectation of a random variable is denoted by $\bb{E}[.]$.  The set
of positive integers is denoted by $\bb{N}$. The ceiling function
and floor function are denoted respectively by $\lc . \rc$ and $\lf
. \rf$ (i.e., $\lc x \rc$ represents the smallest integer no less
than $x$; $\lf x \rf$ represents the largest integer no greater than
$x$).  The gamma function is denoted by $\Ga(.)$.  For any integer
$m$, the combinatoric function $\bi{m}{z}$ with respect to integer
$z$ takes value {\small $\f{ \Ga( m + 1) } { \Ga( z + 1) \Ga (m- z +
1) }$} for $z \leq m$ and value $0$ otherwise. The left limit as
$\ep$ tends to $0$ is denoted as $\lim_{\ep \downarrow 0}$. The
notation ``$\LRA$'' means ``if and only if''.  We use the notation
$\Pr \{ . \mid \se \}$ to indicate that the associated random
samples $X_1, X_2, \cd$ are parameterized by $\se$. The parameter
$\se$ in $\Pr \{ . \mid \se \}$  may be dropped whenever this can be
done without introducing confusion. The other notations will be made
clear as we proceed.

\sect{General Theory}

In this section, we shall develop some general results on the
truncated inverse sampling. Let $X$ be a non-negative random
variable defined in a probability space $(\Om, \mscr{F}, \Pr )$. Our
problem is to estimate the mean, $\mu = \bb{E} [X]$, of $X$ based on
i.i.d. random samples  $X_1, X_2, \cd$ of $X$.   To this end, we
shall adopt a truncated inverse sampling scheme as follows:

Continue sampling until the sample sum is no less than a threshold
value $\ga > 0$ or the number of samples reaches an integer $n$.

Let $\mbf{n}$ be the total number of samples when the sampling is
stopped. By the definition of the truncated inverse sampling scheme,
$\mbf{n}$ is a random variable such that
\[ \mbf{n} (\om) = \min \li \{n, \; \min \{ \ell \in \bb{N} :
\sum_{i = 1}^\ell X_i(\om) \geq \ga \} \ri \}
\]
for any $\om \in \Om$.  Define $\mbf{k} = \sum_{i = 1}^{\mbf{n}}
X_i$.  Then, we can take $\wh{\bs{\mu}} = \f{ \min \{ \mbf{k}, \ga
\} } { \mbf{n} }$ as the estimator for $\mu = \bb{E} [X]$.

With regard to the distribution of $\wh{\bs{\mu}}$, we have \beT For
any $z > 0$,  {\small
\[ \Pr \{ \wh{\bs{\mu}} \leq z \} = \bec \Pr \{ \sum_{i = 1}^{\lc
\ga \sh z \rc - 1} X_i <
\ga \} & \tx{for} \; \ga \leq n z,\\
\Pr \li \{ \sum_{i = 1}^n X_i \leq n z \ri \} & \tx{for} \; \ga
> n z \eec \qqu  \Pr \{  \wh{\bs{\mu}} \geq z \} = \bec \Pr  \{ \sum_{i = 1}^{
\lf \ga \sh z \rf} X_i
\geq \ga  \} & \tx{for} \; \ga \leq n z,\\
\Pr \{ \sum_{i = 1}^{n} X_i \geq n z  \} & \tx{for} \; \ga
> n z. \eec
\]}
\eeT

With regard to the average sample number $\bb{E} [ \mbf{n} ]$ of the
truncated inverse sampling associated with random variable $X$, we
have

\beT For any non-negative random variable $X$ with positive mean and
finite variance, {\small $\bb{E}[ \mbf{n} ] < \min \{ n, \f{\ga}{\mu
} + 1 \}$}.  Specially, if $\ga$ is a positive integer and $X$ is a
Bernoulli random variable such that $\Pr \{ X = 1 \} = 1 - \Pr \{ X
= 0 \} = p \in (0, 1)$, then {\small $\bb{E}[ \mbf{n} ] < \min \{ n,
\f{\ga}{p } \}$}. \eeT

\sect{Estimation of Binomial Parameters}

In this section, we shall consider the estimation of a binomial
parameter based on truncated inverse sampling.  Let $X$ be a
Bernoulli random variable such that $\Pr \{ X = 1 \} = 1 - \Pr \{ X
= 0 \} = p \in (0, 1)$.  Our goal is to estimate $p$ based on i.i.d.
random samples $X_1, X_2, \cd$ of $X$.   Since $X_i$ assumes only
two possible values $0$ or $1$, the threshold value $\ga$ shall be
restricted to an integer.  The estimator for $p$ can be taken as
$\wh{\bs{p}} = \f{ \min \{ \mbf{k}, \ga \} } { \mbf{n} } = \f{
\mbf{k} } { \mbf{n} }$, where $\mbf{k}$ and $\mbf{n}$ have been
defined in Section 2.

In order to estimate $p$ via truncated inverse sampling, a critical
problem is the determination of the threshold value $\ga$ and the
maximum sample size $n$.  By making use of  functions
\[
\mscr{M}_{\mrm{B}}  (z, \mu) =  z \ln \li ( \f{\mu}{z} \ri ) + ( 1 -
z) \ln \li ( \f{ 1 - \mu } { 1 - z } \ri ), \qqu \mscr{M}_{\mrm{I}}
(z, \mu) = \f{1}{z} \mscr{M}_{\mrm{B}} (z, \mu) \] for $0 < z < 1$
and $0 < \mu < 1$, we have derived the following result.

\beT Let $0 < \de < 1$. Let $0 < \vep_a < \vep_r < 1$ be
respectively the margins of absolute and relative errors such that
 $\f{\vep_a}{\vep_r} + \vep_a \leq \f{1}{2}$.  Then,
{\small $\Pr \li \{ | \wh{\bs{p}} - p | < \vep_a \; \mrm{or} \; \li
| \f{\wh{\bs{p}} - p}{p} \ri | < \vep_r  \ri \} > 1 - \de$} provided
that {\small $n > \f{ \ln (\de \sh 2) } { \mscr{M}_{\mrm{B}} (
p^\star + \vep_a, p^\star )  }$} and {\small $\ga > \f{ \ln (\de \sh
2) } { \mscr{M}_{\mrm{I}}  ( p^\star + \vep_a, p^\star ) }$} where
$p^\star = \f{\vep_a}{\vep_r}$. \eeT

Theorem 3 provides explicit formulas for determining the threshold
value $\ga$ and the maximum sample size $n$.   To reduce
conservatism, we can take a computational approach to obtain smaller
$\ga$ and $n$. In this direction, the following theorem which is of
fundamental importance.

 \beT \la{THM1}

Let $\mscr{L}(.)$ and $\mscr{U}(.)$ be monotone functions. Let the
supports of $\mscr{L}(\wh{\bs{p}})$ and $\mscr{U}(\wh{\bs{p}})$ be
denoted by $I_\mscr{L}$ and $I_\mscr{U}$ respectively.  Then, the
maximum of $\Pr \{ p \leq \mscr{L}( \wh{\bs{p}} ) \mid p \}$ with
respect to $p \in [a, b] \subseteq [0, 1]$ is achieved at
$I_\mscr{L} \cap [a, b] \cup \{ a, b\}$ provided that $I_\mscr{L}$
has no closure point in $[a, b]$. Similarly, the maximum of $\Pr \{
p \geq \mscr{U}( \wh{\bs{p}} ) \mid p \}$ with respect to $p \in [a,
b] \subseteq [0, 1]$ is achieved at $I_\mscr{U} \cap [a, b] \cup \{
a, b\}$ provided that $I_\mscr{U}$ has no closure point in $[a, b]$.

\eeT

In Theorem 4, we have used the concept of support. The support of a
random variable is referred to the set of all possible values that
the random variable can assume.  By virtue of Theorem 4, we have
obtained the following results.

 \beT
 Let $0 < \de < 1$ and $\ze > 0$.
 Let $0 < \vep_a < \vep_r < 1$ be respectively the margins of absolute and relative
errors.  Define {\small $n = \li \lf \f{ \ln ( \ze \de) } {
\mscr{M}_{\mrm{B}} ( p^\star + \vep_a, p^\star )  } \ri \rf$} and
{\small $\ga = \li \lf \f{ \ln ( \ze \de) } { \mscr{M}_{\mrm{I}} (
p^\star + \vep_a, p^\star ) } \ri \rf$} with $p^\star =
\f{\vep_a}{\vep_r}$. Define $\mscr{Q}_a^-$ as the support of
$\wh{\bs{p}} - \vep_a$, $\mscr{Q}_a^+$ as the support of
$\wh{\bs{p}} + \vep_a$, $\mscr{Q}_r^+$ as the support of
$\wh{\bs{p}} \sh (1 + \vep_r)$, $\mscr{Q}_r^-$ as the support of
$\wh{\bs{p}} \sh (1 - \vep_r)$. Then, {\small $\Pr \li \{ |
\wh{\bs{p}} - p | < \vep_a \; \mrm{or} \; \li | \f{\wh{\bs{p}} -
p}{p} \ri | < \vep_r \ri \} > 1 - \de$} provided that \bel & & \Pr
\{ \wh{\bs{p}} \geq p + \vep_a \mid p \} \leq \f{\de}{2}, \qqu \fa p
\in
\mscr{Q}_a^-  \cup \{ p^\star \} \cap (0, p^\star] \la{mix1}\\
&   & \Pr \{ \wh{\bs{p}} \leq p - \vep_a \mid p \} \leq \f{\de}{2},
\qqu \fa p \in \mscr{Q}_a^+ \cup \{ p^\star \} \cap (0, p^\star] \la{mix2}\\
 &  &  \Pr \{
\wh{\bs{p}} \geq p (1 + \vep_r) \mid p \} \leq \f{\de}{2}, \qqu \fa
p \in \mscr{Q}_r^+  \cap (p^\star, 1) \la{mix3}\\
&   &  \Pr \{ \wh{\bs{p}} \leq p (1 - \vep_r) \mid p \} \leq
\f{\de}{2}, \qqu \fa p \in \mscr{Q}_r^- \cap (p^\star, 1)
\la{mix4}\eel where these conditions are satisfied when $\ze$ is
smaller than $\f{1}{2}$.
 \eeT

 \bsk

 Clearly, the support of $\wh{\bs{p}}$ is $\{ \f{j}{n} : j = 0, 1,
 \cd, \ga \} \cup \{ \f{\ga}{m} : m = \ga, \ga + 1, \cd, n - 1  \}$.
 Theorem 5 asserts that the prescribed levels of precision and confidence can
 be guaranteed if $\ze$ is small enough.  Hence, we can determine an
 appropriate value of $\ze$ by a bisection search method.

When the sampling is terminated, it is desirable to construct a
confidence interval for $p$.  For this purpose,  we have \beT Let $0
< \de < 1$. Define lower confidence limit $\udl{\bs{p}} \in [0,1)$
such that $\udl{\bs{p}} = 0$ for $\mbf{k} = 0$ and that {\small
$\sum_{i = \mbf{k} }^{ \mbf{n} } \bi{\mbf{n}}{i} \udl{\bs{p}}^i (1 -
\udl{\bs{p}})^{\mbf{n} - i} = \f{\de}{2}$} for $\mbf{k} > 0$. Define
upper confidence limit $\ovl{\bs{p}} \in (0,1]$ such that
$\ovl{\bs{p}} = 1$ for $\mbf{k} = \mbf{n}$ and that {\small $\sum_{i
= 0}^{\mbf{k}} \bi{\mbf{n}}{i} \ovl{\bs{p}}^i (1 -
\ovl{\bs{p}})^{\mbf{n} - i} = \f{\de}{2}$} for $\mbf{k} < \mbf{n}$.
Then, $\Pr \{ \udl{\bs{p}} < p < \ovl{\bs{p}} \} \geq 1 - \de$.

\eeT

\bsk It should be noted the approach of constructing a confidence
interval for $p$ can be considered as a generalization of Clopper
and Pearson's method \cite{Clopper} of interval estimation.

\sect{Estimation of Finite Population Proportion}

In the last section, we have investigated the estimation of a
binomial parameter $p$, which can be considered as the proportion of
an infinite population.  In many situations, the population size is
finite and we shall devote this section to the estimation of the
proportion of a finite population.   Consider a population of $N$
units, among which there are $M$ units having a certain attribute.
It is a frequent problem to estimate the population proportion $p =
\f{M}{N}$ by sampling without replacement. The procedure of sampling
without replacement can be precisely described as follows:

 Each time a single unit is drawn without replacement from the remaining population so
that every unit of the remaining population has equal chance of
being selected.

Such a sampling process can be exactly characterized by random
variables $X_1, \cd, X_N$ defined in a probability space $(\Omega,
\mscr{F}, \Pr)$ such that $X_i$ denotes the characteristics of the
$i$-th sample in the sense that $X_i = 1$ if the $i$-th sample has
the attribute and $X_i = 0$ otherwise.  By the nature of the
sampling procedure, it can be shown that \be \la{dep}
 \Pr \{ X_i =
x_i, \; i = 1, \cd, n \} = \bi{M}{\sum_{i = 1}^n x_i} \bi{N - M}{n -
\sum_{i = 1}^n x_i}
 \li \slash \li [ \bi{n}{\sum_{i = 1}^n x_i} \bi{N}{n} \ri. \ri ]
\ee for any $n \in \{1, \cd, N\}$ and any $x_i \in \{0, 1\}, \; i =
1, \cd, n$.  Moreover, if the proportion $p = \f{M}{N}$ is fixed and
the population size $N$ tends to infinity, the sequence $X_1, X_2,
\cd, X_N$ tends to the i.i.d. random samples of a Bernoulli
variable.

To estimate the population proportion $p$, we can use a sampling
scheme defined by positive integers $\ga$ and $n$ as follows:

Continue sampling without replacement until $\ga$ units found to
have a certain attribute or the number of samples reaches $n$.

Despite the lack of independence in the sequence $X_1, X_2, \cd ,
X_N$ with joint distribution (\ref{dep}), such a sampling method is
also referred to as truncated inverse sampling due to the fact that,
when the sampling is terminated, the number of units having a
certain attribute, denoted by $\mbf{k}$, is actually equal to
$\sum_{i = 1}^{\mbf{n}} X_i$, where $\mbf{n}$ is the sample size
when the sampling is terminated.  This implies that, by relaxing the
independency assumption, we can put such a sampling scheme in the
general framework of truncated inverse sampling described in Section
2. It can be seen that, as the sample size tends to infinity while
the proportion $p$ is being fixed, such a sampling scheme reduces to
the truncated inverse sampling for the estimation of a binomial
parameter as discussed in Section 3.

As in the case of estimating a binomial parameter in Section 3, the
estimator for the proportion of a finite population can be taken as
$\wh{ \bs{p} } = \f{ \min \{ \mbf{k}, \ga \} } { \mbf{n} } = \f{
\mbf{k} } { \mbf{n} }$. In order to determine $n$ and $\ga$ to
guarantee prescribed levels of precision and confidence, we have the
following result.

\beT Let $0 < \de < 1$.  Let $0 < \vep_a < \vep_r < 1$ be
respectively the margins of absolute and relative errors such that
$\f{\vep_a}{\vep_r} + \vep_a \leq \f{1}{2}$.   Then, {\small $\Pr
\li \{ | \wh{\bs{p}} - p | < \vep_a \; \mrm{or} \; \li | \wh{\bs{p}}
- p \ri | < p \vep_r \ri \} > 1 - \de$} provided that {\small $n >
\f{ \ln (\de \sh 2) } { \mscr{M}_{\mrm{B}} ( p^\star + \vep_a,
p^\star )  }$} and {\small $\ga > \f{ \ln (\de \sh 2) } {
\mscr{M}_{\mrm{I}}  ( p^\star + \vep_a, p^\star ) }$}, where
$p^\star = \f{\vep_a}{ \vep_r}$.
 \eeT

Theorem 7 provides explicit formulas for determining threshold value
$\ga$ and the maximum sample size $n$.   To reduce conservatism, we
can take a computational approach to obtain smaller $\ga$ and $n$.
In this direction, the following theorem is useful.

\beT \la{THM3}  Let $\mscr{L}(.)$ and $\mscr{U}(.)$ be
non-decreasing integer-valued functions. Let the supports of
$\mscr{L}(\wh{\bs{p}})$ and $\mscr{U}(\wh{\bs{p}})$ be denoted by
$I_\mscr{L}$ and $I_\mscr{U}$ respectively. Then, the maximum of
$\Pr \{ M \leq \mscr{L}( \wh{\bs{p}} ) \mid M \}$ with respect to $M
\in [a, b] \subseteq [0, N]$, where $a$ and $b$ are integers, is
achieved at $I_\mscr{L} \cap [a, b] \cup \{ a, b\}$. Similarly, the
maximum of $\Pr \{ M \geq \mscr{U}( \wh{\bs{p}} ) \mid M \}$ with
respect to $M \in [a, b]$ is achieved at $I_\mscr{U} \cap [a, b]
\cup \{ a, b\}$. \eeT

By virtue of Theorem 8, we have obtained the following results.

 \beT
 Let $0 < \de < 1$ and $\ze > 0$.  Let $0 < \vep_a < \vep_r < 1$
 be respectively the margins of absolute and relative
errors.  Define {\small $n = \li \lf \f{ \ln ( \ze \de) } {
\mscr{M}_{\mrm{B}} ( p^\star + \vep_a, p^\star )  } \ri \rf$} and
{\small $\ga = \li \lf \f{ \ln ( \ze \de) } { \mscr{M}_{\mrm{I}} (
p^\star + \vep_a, p^\star ) } \ri \rf$} with $p^\star =
\f{\vep_a}{\vep_r}$.  Define $\mscr{Q}_a^-$ as the support of $ \lf
N (\wh{\bs{p}} - \vep_a) \rf$, $\mscr{Q}_a^+$ as the support of $
\lc N (\wh{\bs{p}} + \vep_a) \rc$, $\mscr{Q}_r^+$ as the support of
$ \lf N \wh{\bs{p}} \sh (1 + \vep_r) \rf$, $\mscr{Q}_r^-$ as the
support of $ \lc N \wh{\bs{p}} \sh (1 - \vep_r) \rc$.   Then,
{\small $\Pr \li \{ | \wh{\bs{p}} - p | < \vep_a \; \mrm{or} \; \li
| \wh{\bs{p}} - p \ri | < p \vep_r  \ri \} > 1 - \de$} provided that
\bel & &  \Pr \{ \wh{\bs{p}} \geq p + \vep_a \mid M \} \leq
\f{\de}{2}, \qqu \fa M \in
\mscr{Q}_a^-  \cup \{  \lf N p^\star \rf \} \cap (0, N p^\star] \la{Fmix1}\\
&   & \Pr \{ \wh{\bs{p}} \leq p - \vep_a \mid M \} \leq \f{\de}{2},
\qqu \fa M \in \mscr{Q}_a^+ \cup \{ \lf N p^\star \rf \} \cap (0, N p^\star] \la{Fmix2}\\
 &  &  \Pr \{
\wh{\bs{p}} \geq p (1 + \vep_r) \mid M \} \leq \f{\de}{2}, \qqu \fa
M \in
\mscr{Q}_r^+ \cup \{ \lf N p^\star \rf + 1 \} \cap (N p^\star, N) \la{Fmix3}\\
&   &  \Pr \{ \wh{\bs{p}} \leq p (1 - \vep_r) \mid M \} \leq
\f{\de}{2}, \qqu \fa M \in \mscr{Q}_r^- \cup \{  \lf N p^\star \rf +
1 \} \cap (N p^\star, N) \la{Fmix4}\eel where these conditions are
satisfied when $\ze$ is smaller than $\f{1}{2}$.
 \eeT

 \bsk

 Clearly, the support of $\wh{\bs{p}}$ is $\{ \f{j}{n} : j = 0, 1,
 \cd, \ga \} \cup \{ \f{\ga}{m} : m = \ga, \ga + 1, \cd, n - 1  \}$.
 It is asserted by Theorem 9 that the prescribed levels of precision and confidence can
 be guaranteed if $\ze$ is small enough.  Therefore, an
 appropriate value of $\ze$ can be determined by a bisection search method.

In order to construct a confidence interval for $M$, we have

\beT Let $\bs{M}_l$ be the smallest integer such that {\small
$\sum_{i = \mbf{k} }^{\mbf{n}} \bi{\bs{M}_l}{i} \bi{N -
\bs{M}_l}{\mbf{n} - i} \sh \bi{N}{\mbf{n}}
> \f{\de}{2}$}.  Let $\bs{M}_u$ be the largest
integer such that {\small $\sum_{i = 0}^{\mbf{k}} \bi{\bs{M}_u}{i}
\bi{N - \bs{M}_u}{\mbf{n} - i} \sh \bi{N}{\mbf{n}}
> \f{\de}{2}$}.  Then, $\Pr \{ \bs{M}_l \leq M \leq \bs{M}_u \} \geq
1 - \de$.
 \eeT

With regard to the average sample number $\bb{E} [ \mbf{n} ]$,  we
have

\beT  If the population proportion $p$ is positive, then {\small
$\bb{E}[ \mbf{n} ] < \min \li \{ n, \f{\ga}{p } \ri \}$}. \eeT

\sect{Estimation of Poisson Parameters}

Let $X$ be a Poisson random variable with mean $\lm > 0$. It is a
frequent problem to estimate $\lm$ based on i.i.d. random samples
$X_1, X_2, \cd$ of $X$. This can be accomplished by using the
truncated inverse sampling scheme described in Section 2. Since
$X_i$ is an integer-valued random variable, we shall restrict the
threshold $\ga$ to be a positive integer.   We take $\wh{\bs{\lm}} =
\f{ \min \{ \mbf{k}, \ga \}  } { \mbf{n} }$ as the estimator for
$\lm$, where $\mbf{k}$ and $\mbf{n}$ have been defined in Section 2.

To determine the threshold $\ga$ and the maximum sample size $n$, we
need to have an upper bound for $\lm$.  We do not pursue results
along this line.  We are more interested in the construction of a
confidence interval when the sampling is completed.  For this
purpose, we have

\beT Let $0 < \de < 1$. Define lower confidence limit
$\udl{\bs{\lm}}$ such that $\udl{\bs{\lm}} = 0$ for $\wh{\bs{\lm}} =
0$,  $\sum_{i = \ga}^{\iy} \f{1}{ i! } ( \mbf{n} \udl{\bs{\lm}} )^i
\exp ( - \mbf{n} \udl{\bs{\lm}} ) = \f{\de}{2}$ for $\wh{\bs{\lm}}
\geq \f{\ga}{n}$, and $\sum_{ i = \mbf{k}}^\iy \f{1}{i!} (n
\udl{\bs{\lm}} )^i \exp ( - n \udl{\bs{\lm}} ) = \f{\de}{2}$ for $0
< \wh{\bs{\lm}} < \f{\ga}{n}$, where $\mbf{k} = \sum_{i=1}^n X_i$.
Define upper confidence limit $\ovl{\bs{\lm}}$ such that
$\ovl{\bs{\lm}} = \iy$ for $\mbf{n} = 1$,  {\small $\sum_{i =
0}^{\ga - 1} \f{1}{i!} [( \mbf{n} - 1) \ovl{\bs{\lm}} ]^i \exp ( - (
\mbf{n} - 1) \ovl{\bs{\lm}} ) = \f{\de}{2}$} for $\f{\ga}{n} \leq
\wh{\bs{\lm}} < \ga$, and $\sum_{ i = 0 }^{\mbf{k}} \f{1}{i!} (n
\ovl{\bs{\lm}} )^i \exp ( - n \ovl{\bs{\lm}} ) = \f{\de}{2}$ for
$\wh{\bs{\lm}} < \f{\ga}{n}$.  Then, $\Pr \{ \udl{\bs{\lm}} < \lm <
\ovl{\bs{\lm}} \} \geq 1 - \de$.

\eeT

\bsk

It should be noted that the interval estimation method described in
Theorem 12 is a generalization of Garwood's interval estimation
method \cite{Garwood}.

\sect{Estimation of Bounded-Variable Means}

Let $X$ be a random variable bounded in $[0, 1]$ with mean $\mu =
\bb{E} [ X ]$.  In many situations, it is desirable to estimate
$\mu$ based on i.i.d. random samples $X_1, X_2, \cd$ of $X$ (see,
e.g., \cite{Dagum} and the references therein). To fulfill this
goal, we shall make use of the truncated inverse sampling scheme
described in Section 2. In order to determine the threshold $\ga$
and maximum sample size $n$ to guarantee prescribed levels of
precision and confidence, we have

\beT Let $0 < \de < 1$. Let $0 < \vep_a < \vep_r < 1$ be
respectively the margins of absolute and relative errors such that
$p^\star + \vep_a \leq \f{1}{2}$ with $p^\star =
\f{\vep_a}{\vep_r}$.  Then, {\small $\Pr \li \{ | \wh{\bs{\mu}} -
\mu | < \vep_a \; \mrm{or} \; \li | \f{\wh{\bs{\mu}} - \mu}{\mu} \ri
| < \vep_r  \ri \} > 1 - \de$} provided that
\[
\ga > \f{1 - \vep_r}{\vep_r}, \qu \ga > \f{ \ln \f{\de}{2} } {
\mscr{M}_{\mrm{I}}  \li ( \f{\ga (p^\star - \vep_a)}{\ga - 1 +
\vep_r}, p^\star \ri )  }, \qu \ga > \f{ \ln \f{\de}{2} } {
\mscr{M}_{\mrm{I}}  ( p^\star + \vep_a, p^\star ) }, \qu n > \f{ \ln
\f{\de}{2} } { \mscr{M}_{\mrm{B}} ( p^\star + \vep_a, p^\star )  }.
\]
\eeT

\bsk

With regard to the interval estimation of $\mu$, we have

\beT

Let $0 < \de < 1$. Define lower confidence limit $\udl{\bs{\mu}} \in
[0, \wh{\bs{\mu}}]$ such that $\udl{\bs{\mu}} = 0$ for
$\wh{\bs{\mu}} = 0$, $\mscr{M}_{\mrm{B}} ( \wh{\bs{\mu}},
\udl{\bs{\mu}}) = \f{ \ln \f{\de}{2} }{n}$ for $0 < \wh{\bs{\mu}} <
\f{\ga}{n}$, and $\mscr{M}_{\mrm{I}}  ( \wh{\bs{\mu}},
\udl{\bs{\mu}})  = \f{ \ln \f{\de}{2} }{\ga}$ for $\wh{\bs{\mu}}
\geq \f{\ga}{n}$.  Define upper confidence limit $\ovl{\bs{\mu}} \in
[ \wh{\bs{\mu}}, 1]$ and that $\ovl{\bs{\mu}} = 1$ for
$\wh{\bs{\mu}} \geq \f{\ga}{\ga + 1}$, $\mscr{M}_{\mrm{B}} (
\wh{\bs{\mu}}, \ovl{\bs{\mu}}) = \f{ \ln \f{\de}{2} }{n}$ for
$\wh{\bs{\mu}} < \f{\ga}{n}$, $\mscr{M}_{\mrm{I}} \li (
\f{\wh{\bs{\mu}} \ga}{ \ga  - \wh{\bs{\mu}} } , \ovl{\bs{\mu}} \ri )
= \f{ \ln \f{\de}{2} }{\ga}$ and $ \f{\wh{\bs{\mu}} \ga}{ \ga  -
\wh{\bs{\mu}} } < \ovl{\bs{\mu}} < 1 $ for $\f{\ga}{\ga + 1} >
\wh{\bs{\mu}} \geq \f{\ga}{n}$. Then, $\Pr \{ \udl{\bs{\mu}} < \mu <
\ovl{\bs{\mu}} \} \geq 1 - \de$.

\eeT

\sect{Conclusion}

In this paper, we have established a general theory of truncated
inverse sampling for estimating the mean value of a large class of
random variables. We have applied such a theory to the common
important variables such as binomial, Poisson, hyper-geometrical,
and bounded variables. Rigorous methods have been derived for
determining the thresholds and maximum sample sizes to ensure
statistical accuracy. Interval estimation methods have also been
developed.

\bsk

\appendix

\sect{Proof of Theorem 1}

The theorem can be shown by establishing Lemmas 1 to 4 as follows.

\beL Suppose $z \geq \f{\ga}{n}$.  Then, $\Pr \{  \wh{\bs{\mu}} \leq
z \} = \Pr \{ \sum_{i = 1}^m  X_i < \ga \}$ where $m = \li \lc
\f{\ga}{z} \ri \rc - 1$.  \eeL

\bpf

By the assumption that $z \geq \f{\ga}{n}$ and the definition of the
sampling scheme, {\small \bee \li \{ \wh{\bs{\mu}}
> z, \; \mbf{k} < \ga \ri \}  =  \li \{  \f{
\mbf{k} } { \mbf{n} }  > z, \; \mbf{k} < \ga \ri \} \subseteq  \li
\{  \f{ \mbf{k} } { \mbf{n} }  > \f{\ga}{n}, \; \mbf{k} < \ga \ri \}
=  \li \{  \f{ \mbf{k} } { \mbf{n} }  > \f{\ga}{n}, \; \mbf{k} <
\ga, \; \mbf{n} = n \ri \} =  \emptyset. \eee} Therefore, {\small
$\{  \wh{\bs{\mu}} > z \}  =  \{ \wh{\bs{\mu}}
> z, \; \mbf{k} < \ga \} \cup \{
\wh{\bs{\mu}} > z, \; \mbf{k} \geq \ga \} = \{ \wh{\bs{\mu}} > z, \;
\mbf{k} \geq \ga \} = \{ \f{ \ga  } { \mbf{n} } > z, \; \mbf{k} \geq
\ga \}$}.  To show the lemma, it remains to show that {\small $\li
\{ \f{ \ga } {  \mbf{n} } > z, \; \mbf{k} \geq \ga \ri \} = \Pr \{
\sum_{i=1}^m X_i \geq \ga \}$}. Since all $X_i$ are non-negative, we
have {\small $\{  \f{ \ga } { \mbf{n} } > z, \; \mbf{k} \geq \ga \}
=  \{ \mbf{n} \leq m, \; \mbf{k} \geq \ga \} \subseteq \{ \mbf{n}
\leq m, \; \sum_{i = 1}^{m} X_i \geq \ga \} \subseteq \{ \sum_{i =
1}^{m} X_i \geq \ga \}$}.  On the other hand, by the assumption that
$z \geq \f{\ga}{n}$, we have $m = \li \lc \f{\ga}{z} \ri \rc - 1
\leq n - 1$.  Hence, by the definition of the sampling scheme, we
have {\small $\{ \sum_{i = 1}^{m} X_i \geq \ga \} \subseteq \{
\mbf{n} \leq m, \; \mbf{k} \geq \ga, \; \sum_{i = 1}^{m} X_i \geq
\ga \} \subseteq  \{ \mbf{n} \leq m, \; \mbf{k} \geq \ga \}$}. It
follows that {\small $\{ \sum_{i = 1}^{m} X_i \geq \ga \} = \{
\mbf{n} \leq m, \; \mbf{k} \geq \ga \} = \{  \f{ \ga  } { \mbf{n} }
> z, \; \mbf{k} \geq \ga \} = \{ \wh{\bs{\mu}}
> z \}$}.  This completes the proof of the lemma.

\epf

\beL  Suppose $z < \f{\ga}{n}$.  Then, $\Pr \li \{  \wh{\bs{\mu}}
\leq z \ri \} =  \Pr \li \{ \f{ \sum_{i = 1}^n X_i}{n} \leq z \ri
\}$.
 \eeL

 \bpf  By the assumption that $z < \f{\ga}{n}$ and the definition of the
sampling scheme, \bee \li \{ \wh{\bs{\mu}} \leq z, \; \mbf{k} \geq
\ga \ri \}  = \li \{  \f{ \ga } { \mbf{n} } \leq z, \; \mbf{k} \geq
\ga \ri \} \subseteq \li \{  \f{ \ga  } { \mbf{n} } < \f{\ga}{n}, \;
\mbf{k} \geq \ga \ri \} =  \li \{ \mbf{n}  > n, \; \mbf{k} \geq \ga
\ri \} = \emptyset. \eee Therefore,  {\small \bee \{ \wh{\bs{\mu}}
\leq z \} & = & \li \{ \wh{\bs{\mu}} \leq z, \; \mbf{k} < \ga \ri \}
\cup \li \{ \wh{\bs{\mu}} \leq z, \; \mbf{k} \geq \ga \ri \} =  \li
\{ \wh{\bs{\mu}} \leq z, \; \mbf{k} < \ga \ri \} = \li \{ \f{
\mbf{k} } { \mbf{n} } \leq z, \; \mbf{k} < \ga \ri \}\\
&  = & \li \{ \f{ \mbf{k} } { \mbf{n} } \leq z, \; \mbf{k} < \ga, \;
\mbf{n} = n \ri \}  =  \li \{  \f{ \sum_{i = 1}^{n} X_i  } { n }
\leq z, \; \sum_{i = 1}^{n} X_i < \ga, \; \mbf{n} = n \ri \}
\subseteq \li \{  \f{ \sum_{i = 1}^{n} X_i  } { n } \leq z \ri \}.
\eee} On the other hand, by the definition of the sampling scheme
and the assumption that $z < \f{\ga}{n}$, we have {\small $\li \{
\f{ \sum_{i = 1}^{n} X_i } { n } \leq z \ri \}  \subseteq  \li \{
\f{ \sum_{i = 1}^{n} X_i } { n } \leq z, \; \; \mbf{n} = n \ri \} =
\li \{  \f{ \sum_{i = 1}^{n} X_i  } { n } \leq z, \; \sum_{i =
1}^{n} X_i < \ga, \; \mbf{n} = n \ri \} =  \{  \wh{\bs{\mu}} \leq z
\}$}. It follows that {\small $\li \{  \wh{\bs{\mu}} \leq z \ri \} =
\li \{ \f{\sum_{i = 1}^{n} X_i}{n} \leq z \ri \}$}.  This completes
the proof of the lemma.
 \epf

\beL  Suppose $z \geq \f{\ga}{n}$.  Then, $\Pr \{  \wh{\bs{\mu}}
\geq z \} = \Pr \{ \sum_{i = 1}^{m} X_i \geq \ga \}$ where $m = \li
\lf \f{\ga}{z} \ri \rf$.
 \eeL

 \bpf
By the assumption that $z \geq \f{\ga}{n}$ and the definition of the
sampling scheme, {\small \bee \li \{ \wh{\bs{\mu}} \geq z, \;
\mbf{k} < \ga \ri \}  =  \li \{  \f{ \mbf{k} } { \mbf{n} } \geq z,
\; \mbf{k} < \ga \ri \} = \li \{  \f{ \mbf{k} } { \mbf{n} } \geq z,
\; \mbf{k} < \ga, \; \mbf{n} = n \ri \}  = \emptyset.  \eee}
Therefore,  {\small $ \{ \wh{\bs{\mu}} \geq z \}  =  \{
\wh{\bs{\mu}} \geq z, \; \mbf{k} < \ga \} \cup \{ \wh{\bs{\mu}} \geq
z, \; \mbf{k} \geq \ga \} = \{ \wh{\bs{\mu}} \geq z, \; \mbf{k} \geq
\ga \} = \{ \f{ \ga  } { \mbf{n} } \geq z, \; \mbf{k} \geq \ga \}$}.
To show the lemma, it remains to show that {\small $\{ \f{ \ga  } {
\mbf{n} } \geq z, \; \mbf{k} \geq \ga \} = \Pr \{ \sum_{i = 1}^m X_i
\geq \ga \}$}. Since all $X_i$ are non-negative, we have {\small $\{
\f{ \ga } { \mbf{n} } \geq z, \; \mbf{k} \geq \ga \} = \{ \mbf{n}
\leq m, \; \mbf{k} \geq \ga \} \subseteq \{ \mbf{n} \leq m, \;
\sum_{i = 1}^{m} X_i \geq \ga \} \subseteq \{ \sum_{i = 1}^{m} X_i
\geq \ga \}$}.  On the other hand, by the assumption that $z \geq
\f{\ga}{n}$, we have $m = \li \lf \f{\ga}{z} \ri \rf \leq n$. By the
definition of the sampling scheme and the fact that all $X_i$ are
non-negative, {\small $ \{ \sum_{i = 1}^{m} X_i \geq \ga \}
\subseteq \{ \mbf{n} \leq m, \; \sum_{i = 1}^m X_i \geq \ga \}
\subseteq \{ \mbf{n} \leq m, \; \mbf{k} \geq \ga \}$}. Hence,
{\small $\{  \sum_{i = 1}^{m} X_i \geq \ga \} = \{ \mbf{n} \leq m,
\; \mbf{k} \geq \ga \} = \{ \f{ \ga } { \mbf{n} } \geq z, \; \mbf{k}
\geq \ga \} = \{ \wh{\bs{\mu}} \geq z \}$}. This completes the proof
of the lemma.

 \epf

\beL Suppose $z < \f{\ga}{n}$.  Then, $\Pr \{  \wh{\bs{\mu}} \geq z
\} = \Pr \li \{  \f{ \sum_{i = 1}^{n} X_i  } { n } \geq z \ri \}$.
\eeL

\bpf By the assumption that $z < \f{\ga}{n}$ and the definition of
the sampling scheme, \bee \li \{ \wh{\bs{\mu}} < z, \; \mbf{k} \geq
\ga \ri \}  =  \li \{  \f{ \ga } { \mbf{n} } < z, \; \mbf{k} \geq
\ga \ri \} \subseteq \li \{ \f{ \ga  } { \mbf{n} } < \f{\ga}{n}, \;
\mbf{k} \geq \ga \ri \}  =  \li \{ \mbf{n}  > n, \; \mbf{k} \geq \ga
\ri \} =  \emptyset. \eee  Therefore, {\small \bee \{ \wh{\bs{\mu}}
< z \} & = & \li \{ \wh{\bs{\mu}} < z, \; \mbf{k} < \ga \ri \} \cup
\li \{ \wh{\bs{\mu}} < z, \; \mbf{k} \geq \ga \ri \}  =  \li \{
\wh{\bs{\mu}} < z, \; \mbf{k} < \ga \ri \} = \li \{ \f{ \mbf{k} } {
\mbf{n} } < z, \; \mbf{k} < \ga \ri \}\\
& = & \li \{  \f{ \mbf{k}  } { \mbf{n} } < z, \; \mbf{k} < \ga, \;
\mbf{n} = n \ri \} = \li \{  \f{ \sum_{i = 1}^{n} X_i  } { n } < z,
\; \sum_{i = 1}^{n} X_i < \ga, \; \mbf{n} = n \ri \} \subseteq \li
\{  \f{ \sum_{i = 1}^{n} X_i  } { n } < z \ri \}. \eee} On the other
hand, by the definition of the sampling scheme and the assumption
that $z < \f{\ga}{n}$, {\small \bee \li \{  \f{ \sum_{i = 1}^{n} X_i
} { n } < z \ri \}  \subseteq \li \{  \f{ \sum_{i = 1}^{n} X_i  } {
n } < z, \; \; \mbf{n} = n \ri \} =  \li \{  \f{ \sum_{i = 1}^{n}
X_i } { n } < z, \; \sum_{i = 1}^{n} X_i < \ga, \; \mbf{n} = n \ri
\} = \{ \wh{\bs{\mu}} < z \}.  \eee} It follows that $\li \{
\wh{\bs{\mu}} < z \ri \} =  \li \{ \f{\sum_{i = 1}^{n} X_i}{n} < z
\ri \}$, i.e., $\li \{ \wh{\bs{\mu}} \geq z \ri \} = \li \{
\f{\sum_{i = 1}^{n} X_i}{n} \geq z \ri \}$. This completes the proof
of the lemma.

\epf

\sect{Proof of Theorem 2}

By the definition of the truncated inverse sampling scheme, \bee
\bb{E} [ \mbf{n} ] & = & n \Pr \li \{ \sum_{i = 1}^n X_i < \ga \ri
\}  + \sum_{m = 1}^n m \Pr \li \{ \sum_{i = 1}^{m - 1} X_i <
\ga, \; \sum_{i = 1}^m X_i \geq \ga \ri \}\\
& < & n \Pr \li \{ \sum_{i = 1}^n X_i < \ga \ri \}  + \sum_{m = 1}^n
n \Pr \li \{ \sum_{i = 1}^{m - 1} X_i < \ga, \; \sum_{i = 1}^m X_i
\geq \ga \ri \} = n.  \eee By the fact that  $X_i$ is non-negative,
{\small \[ \li ( \bigcup_{m = n + 1}^\iy \li \{ \sum_{i = 1}^{m - 1}
X_i < \ga, \; \sum_{i = 1}^m X_i \geq \ga \ri \}  \ri ) \bigcup \li
\{ \sum_{i = 1}^\iy X_i < \ga \ri \} = \li \{ \sum_{i = 1}^n X_i <
\ga \ri \}.
\]}
Since $\bb{E} [ X ] = \mu$ is positive and the corresponding
variance $\si^2$ is finite, we have,  by Chebyshev's inequality,
{\small \bee 0 \leq \Pr \li \{ \sum_{i = 1}^\iy X_i < \ga \ri \} & =
& \lim_{k \to \iy} \Pr \li \{ \sum_{i = 1}^k X_i < \ga \ri \} =
\lim_{k \to \iy} \Pr \li \{ \f{ \sum_{i = 1}^k X_i}{k} - \mu <
\f{\ga}{k} - \mu \ri \}\\
& \leq & \lim_{k \to \iy} \Pr \li \{ \li |  \f{ \sum_{i = 1}^k
X_i}{k} - \mu \ri | > \li | \f{\ga}{k} - \mu \ri | \ri \} \leq
\lim_{k \to \iy} \f{ \f{\si^2}{k} } { \li | \f{\ga}{k} - \mu \ri |^2
} = 0.
 \eee}
Hence, {\small \bee \bb{E} [ \mbf{n} ] & = &  \sum_{m = n + 1}^\iy n
\Pr \li \{ \sum_{i = 1}^{m - 1} X_i < \ga, \; \sum_{i = 1}^m X_i
\geq \ga \ri \}   + \sum_{m = 1}^n m \Pr \li \{ \sum_{i = 1}^{m - 1}
X_i < \ga, \; \sum_{i = 1}^m X_i \geq \ga \ri \}\\
& < & \sum_{m = 1}^\iy m \Pr \li \{ \sum_{i = 1}^{m - 1} X_i < \ga,
\; \sum_{i = 1}^m X_i \geq \ga \ri \} = \bb{E} [ \mbf{m} ], \eee}
where $\mbf{m}$ is the sample number of the classical {\it inverse
sampling} scheme with the following stopping rule: Sampling is
continued until the sample sum is no less than $\ga$.  By the
definition of the classical inverse sampling, we have $\sum_{i =
1}^{\mbf{m} - 1} X_i < \ga$. Applying Wald's equation, we have
$\bb{E} [ \sum_{i = 1}^{\mbf{m} - 1} X_i ] = \bb{E} [ \mbf{m} - 1 ]
\; \bb{E} [ X ] < \ga$, which implies that $\bb{E} [ \mbf{m} ] <
\f{\ga}{\bb{E} [ X ]} + 1 = \f{\ga}{\mu} + 1$.  Since $\bb{E} [
\mbf{n} ]$ is less than both $n$ and $\bb{E} [ \mbf{m} ]$ as shown
above, we have {\small $\bb{E} [ \mbf{n} ] < \min \{ n, \f{\ga}{\mu}
+ 1 \}$}.

In the special case that $\ga$ is a positive integer and that $X$ is
a Bernoulli random variable such that $\bb{E} [ X ] = p \in (0, 1)$,
we have $\sum_{i = 1}^{\mbf{m}} X_i = \ga$ and consequently, by
Wald's equation,  $\bb{E} [ \sum_{i = 1}^{\mbf{m}} X_i ] = \bb{E} [
\mbf{m} ] \; \bb{E} [ X ] = \ga$, from which we get $\bb{E} [
\mbf{m} ] = \f{\ga}{\bb{E} [ X ]} = \f{\ga}{p}$ and it follows that
{\small $\bb{E} [ \mbf{n} ] < \min \{ n, \f{\ga}{p} \}$}.  This
completes the proof the theorem.

\sect{Proof of Theorem 3}

We need some preliminary results.   The following lemma is a slight
modification of Hoeffding \cite{Hoeffding}.

\beL \la{lem1} Let $X_1, \cd, X_n$ be i.i.d. random variables
bounded in $[0,1]$ with common mean value $\mu \in (0,1)$. Then,
{\small $\Pr \li \{ \f{\sum_{i =1}^n X_i}{n} \geq z \ri \} \leq \exp
\li ( n  \mscr{M}_{\mrm{B}} (z, \mu) \ri )$ } for $1 \geq z \geq
\mu$. Similarly, {\small $\Pr \li \{  \f{\sum_{i =1}^n X_i}{n} \leq
z  \ri \} \leq \exp \li (  n \mscr{M}_{\mrm{B}} (z, \mu) \ri )$ }
for $0 \leq z \leq \mu$. \eeL

\bpf For $z = \mu$, we have {\small $\Pr \li \{ \f{\sum_{i =1}^n
X_i}{n} \geq z \ri \} \leq \exp \li ( n \mscr{M}_{\mrm{B}} (z, \mu)
\ri ) = 1$}.  For $\mu < z < 1$, it was shown by Hoeffding in
\cite{Hoeffding} that {\small $\Pr \li \{ \f{\sum_{i =1}^n X_i}{n}
\geq z \ri \} \leq \exp \li ( n \mscr{M}_{\mrm{B}} (z, \mu) \ri )$}.
For $z = 1$, we have {\small $\Pr \li \{ \f{\sum_{i =1}^n X_i}{n}
\geq z \ri \} = \prod_{i = 1}^n \Pr \{ X_i = 1 \} \leq \prod_{i =
1}^n \bb{E} [ X_i ] = \mu^n = \exp \li ( n \mscr{M}_{\mrm{B}} (1,
\mu) \ri )$}.

For $z = 0$, we have {\small $\Pr \li \{ \f{\sum_{i =1}^n X_i}{n}
\leq z \ri \} = \prod_{i = 1}^n (1 - \Pr \{ X_i \neq 0 \} ) \leq
\prod_{i = 1}^n (1 - \bb{E} [ X_i ] ) = (1 - \mu)^n = \exp \li ( n
\mscr{M}_{\mrm{B}} (0, \mu) \ri )$}.  For $0 < z < \mu$, it was
shown by Hoeffding in \cite{Hoeffding} that {\small $\Pr \li \{
\f{\sum_{i =1}^n X_i}{n} \leq z \ri \} \leq \exp \li ( n
\mscr{M}_{\mrm{B}} (z, \mu) \ri )$}. For $z = \mu$, we have {\small
$\Pr \li \{ \f{\sum_{i =1}^n X_i}{n} \leq z \ri \} \leq \exp \li ( n
\mscr{M}_{\mrm{B}} (z, \mu) \ri ) = 1$}. \epf

\beL \la{lem3} Let $0 < \vep < 1$.  Then, $\mscr{M}_{\mrm{I}} ( \mu
+ \vep \mu, \mu)$ is monotonically decreasing with respect to $\mu
\in \li (0, \f{1}{1 + \vep} \ri )$. Similarly, $\mscr{M}_{\mrm{I}}(
\mu - \vep \mu, \mu)$ is monotonically decreasing with respect to
$\mu \in (0, 1)$. \eeL

\bpf

Note that {\small $\f{ \pa \mscr{M}_{\mrm{I}} ( \mu + \vep \mu,
\mu)} { \pa \mu }  = - \f{1}{\mu^2 (1 + \vep)}  \ln \li [ \f{1 -
\mu}{1 - \mu (1 + \vep) } \ri ] + \f{\vep}{\mu (1 - \mu) (1 + \vep)}
\leq 0$ } if {\small $\ln \li [ \f{1 - \mu}{1 - \mu (1 + \vep) } \ri
]  \geq \f{\vep \mu }{1 - \mu}$}, i.e., {\small \be \la{con} \ln \li
(1 - \f{\vep \mu }{1 - \mu} \ri) \leq - \f{\vep \mu }{1 - \mu}. \ee}
As a consequence of $0 < \mu < \f{1}{1 + \vep}$, we have $0 <
\f{\vep \mu }{1 - \mu} < 1$. Since $\ln (1 - x) < - x$ for any $x
\in (0,1)$, it follows that (\ref{con}) holds and thus
$\mscr{M}_{\mrm{I}} ( \mu + \vep \mu, \mu)$ is monotonically
decreasing with respect to {\small $\mu \in \li (0, \f{1}{1 + \vep}
\ri )$}.

Similarly, to show that $\mscr{M}_{\mrm{I}}( \mu - \vep \mu, \mu)$
is monotonically decreasing with respect to $\mu$, note that {\small
$\f{ \pa \mscr{M}_{\mrm{I}}( \mu - \vep \mu, \mu) } { \pa \mu }  = -
\f{1}{\mu^2 (1 - \vep)} \ln \li [ \f{1 - \mu}{1 - \mu (1 - \vep) }
\ri ]  - \f{\vep}{\mu (1 - \mu) (1 - \vep)} \leq 0$} if {\small $\ln
\li [ \f{1 - \mu}{1 - \mu (1 - \vep) } \ri ] \geq - \f{\vep \mu }{1
- \mu}$}, i.e., {\small \be \la{conze} \ln \li (1 + \f{\vep \mu }{1
- \mu} \ri) \leq \f{\vep \mu }{1 - \mu}. \ee}
 Since $\f{\vep \mu }{1 - \mu} > 0$ and $\ln (1 + x) <  x$ for any $x \in (0,\iy)$, we
 have that (\ref{conze}) holds and thus $\mscr{M}_{\mrm{I}}( \mu - \vep \mu,
\mu)$ is monotonically decreasing with respect to $\mu \in (0, 1)$.

\epf

\beL \la{lem4} $\mscr{M}_{\mrm{I}} ( \mu + \vep \mu, \mu) >
\mscr{M}_{\mrm{I}}( \mu - \vep \mu, \mu)$ for $\mu \in \li (0,
\f{1}{2} \ri )$ and $0 < \vep < 1$. \eeL

\bpf

Direct computation shows that {\small \[ \f{ \pa \mscr{M}_{\mrm{I}}
( \mu + \vep \mu, \mu)} { \pa \vep }  =  - \f{1}{(1 + \vep)^2 \mu}
\ln \li [ \f{1 - \mu}{ 1 - (1 + \vep) \mu } \ri ], \qqu  \f{ \pa
\mscr{M}_{\mrm{I}}( \mu - \vep \mu, \mu)} { \pa \vep } =   \f{1}{(1
- \vep)^2 \mu} \ln \li [  \f{1 - \mu}{ 1 - (1 - \vep) \mu } \ri ].
\]}
Since {\small $\ln \li [  \f{1 - \mu}{ 1 - (1 + \vep) \mu } \ri ]  <
\f{ \vep \mu } { 1 - (1 + \vep) \mu  }$} and {\small $\ln \li [ \f{1
- \mu}{ 1 - (1 - \vep) \mu } \ri ] < - \f{ \vep \mu } { 1 - (1 -
\vep) \mu }$}, we have {\small \[ \f{ \pa \mscr{M}_{\mrm{I}} ( \mu +
\vep \mu, \mu)} { \pa \vep } - \f{ \pa \mscr{M}_{\mrm{I}}( \mu -
\vep \mu, \mu)} { \pa \vep }  > - \f{1}{(1 + \vep)^2 \mu} \f{ \vep
\mu } { 1 - (1 + \vep) \mu } + \f{1}{(1 - \vep)^2 \mu} \f{ \vep
\mu}{1 - (1 - \vep) \mu} > 0
\]}
if $(1 + \vep)^2 [1 - (1 + \vep) \mu ] - (1 - \vep)^2 [ 1 - (1 -
\vep) \mu ] > 0$,  or equivalently, $4 \vep - 2 \vep ( 3 + \vep^2 )
\mu > 0$, which is true because $4 \vep - 2 \vep ( 3 + \vep^2 ) \mu
> 4 \vep - 2 \vep ( 3 + 1 ) \times \f{1}{2}  = 0$ as a result of $0
< \vep < 1$ and $0 < \mu < \f{1}{2}$.  The lemma immediately follows
from the fact that $\f{ \pa \mscr{M}_{\mrm{I}} ( \mu + \vep \mu,
\mu)} { \pa \vep }$ is greater than $\f{ \pa \mscr{M}_{\mrm{I}}( \mu
- \vep \mu, \mu)} { \pa \vep }$  for $0 < \vep < 1, \; 0 < \mu <
\f{1}{2}$ and $\mscr{M}_{\mrm{I}} ( \mu + \vep \mu, \mu)$ is equal
to $\mscr{M}_{\mrm{I}}( \mu - \vep \mu, \mu)$ for $\vep = 0$.

\epf

\beL \la{lem5} Let $0 < \vep < \f{1}{2}$. Then,
$\mscr{M}_{\mrm{I}}(\mu + \vep, \mu)$ is monotonically increasing
with respect to $\mu \in \li (0, \f{1}{2} - \vep \ri )$. Similarly,
$\mscr{M}_{\mrm{I}}(\mu - \vep, \mu)$ is monotonically increasing
with respect to $\mu \in \li (\vep, \f{1}{2} + \vep \ri )$. \eeL

\bpf

It can be shown that {\small $\f{ \mscr{M}_{\mrm{I}}(\mu + \vep,
\mu) }{\pa \mu}  = - \f{ 1 } { (\mu + \vep)^2 } \ln \li ( \f{1 - \mu
} { 1 - \mu - \vep } \ri ) + \f{ \vep } { \mu(\mu + \vep)(1 - \mu) }
> 0$} if {\small $\ln \li ( \f{1 - \mu} { 1 - \mu - \vep  } \ri ) < \f{ \vep (\mu + \vep) }
{ \mu(1 - \mu) }$}.  Since {\small $\ln \li ( \f{1 - \mu } { 1 - \mu
- \vep } \ri ) < \f{\vep}{ 1 - \mu - \vep }$},  it suffices to have
{\small $\f{\vep}{ 1 - \mu - \vep } < \f{ \vep (\mu + \vep) } {
\mu(1 - \mu) }$}, or equivalently, $\mu(1 - \mu)  < (1 - \mu - \vep)
(\mu + \vep) = \mu(1 - \mu) - \vep \mu + (1 - \mu) \vep - \vep^2$,
which can be ensured by $0< \mu < \f{1}{2} - \vep$.  This proves the
first statement of the lemma.

Similarly, {\small $\f{ \mscr{M}_{\mrm{I}}(\mu - \vep, \mu) }{\pa
\mu} = - \f{ 1 } { (\mu - \vep)^2 } \ln \li ( \f{1 - \mu } { 1 - \mu
+ \vep } \ri ) - \f{ \vep } { \mu(\mu - \vep)(1 - \mu) } > 0$} if
{\small $\ln \li ( \f{1 - \mu} { 1 - \mu + \vep } \ri ) < - \f{ \vep
(\mu - \vep) } { \mu(1 - \mu) }$}. Since {\small $\ln \li ( \f{1 -
\mu } { 1 - \mu + \vep } \ri ) < - \f{ \vep }{1 - \mu + \vep }$}, to
ensure {\small $\f{ \mscr{M}_{\mrm{I}}(\mu - \vep, \mu) }{\pa \mu} >
0$}, it suffices to have {\small $- \f{ \vep }{1 - \mu + \vep } < -
\f{ \vep (\mu - \vep) } { \mu(1 - \mu) }$}, or equivalently, $\mu(1
- \mu) > (1 - \mu + \vep) (\mu - \vep) = \mu(1 - \mu) + \vep \mu -
\vep (1 - \mu) - \vep^2$, which can be guaranteed by $\vep < \mu <
\f{1}{2} + \f{\vep}{2}$. This proves the second statement of the
lemma.

\epf

\beL \la{lem6} $\mscr{M}_{\mrm{I}}(\mu + \vep, \mu)
> \mscr{M}_{\mrm{I}}(\mu - \vep, \mu)$  for  $0 < \vep < \mu <
\f{1}{2}$. \eeL

\bpf

It can be verified that {\small $\f{ \pa \mscr{M}_{\mrm{I}}(\mu +
\vep, \mu) } { \pa \vep }
 =  - \f{ 1 } { (\mu + \vep)^2 }  \ln \li ( \f{1 - \mu } { 1 - \mu - \vep
} \ri )$}  and {\small $\f{ \pa \mscr{M}_{\mrm{I}}(\mu - \vep, \mu)
} { \pa \vep } = \f{ 1 } { (\mu - \vep)^2 }  \ln \li ( \f{1 - \mu }
{ 1 - \mu + \vep } \ri )$}.  Since {\small $\ln \li ( \f{1 - \mu } {
1 - \mu - \vep  } \ri ) < \f{\vep}{ 1 - \mu - \vep  }$} and {\small
$\ln \li ( \f{1 - \mu } { 1 - \mu + \vep } \ri ) < - \f{ \vep }{1 -
\mu + \vep }$}, to ensure {\small $\f{ \pa \mscr{M}_{\mrm{I}}(\mu +
\vep, \mu) } { \pa \vep } > \f{ \pa \mscr{M}_{\mrm{I}}(\mu - \vep,
\mu) } { \pa \vep }$}, it suffices to have {\small $- \f{ 1 } { (\mu
+ \vep)^2 } \f{\vep}{ 1 - (\mu + \vep)  }
> - \f{ 1 } { (\mu - \vep)^2 }  \f{ \vep }{1 - \mu + \vep }$}, or equivalently, $2 \mu - 3 \mu^2
> \vep^2$, which is true because
$2 \mu - 3 \mu^2 - \vep^2 > 2 \mu - 3 \mu^2 - \mu^2 = 2 \mu ( 1 - 2 \mu) > 0$
 as a result of $0 < \vep < \mu < \f{1}{2}$.  Therefore, the
lemma is true since $\mscr{M}_{\mrm{I}}(\mu + \vep, \mu) =
\mscr{M}_{\mrm{I}}(\mu - \vep, \mu)$ for $\vep = 0$ and {\small $\f{
\pa \mscr{M}_{\mrm{I}}(\mu + \vep, \mu) } {  \pa \vep } > \f{ \pa
\mscr{M}_{\mrm{I}}(\mu - \vep, \mu) } {  \pa \vep }$ for $0 < \vep <
\mu < \f{1}{2}$}.  This completes the proof of the lemma.

\epf

\beL \la{lem7} Let $0 < \vep < \f{1}{2}$. Then, $\mscr{M}_{\mrm{B}}
(\mu + \vep, \mu)$ is monotonically increasing with respect to $\mu
\in \li (0,  \f{1}{2} - \vep \ri )$. Similarly, $\mscr{M}_{\mrm{B}}
(\mu - \vep, \mu)$ is monotonically increasing with respect to $\mu
\in \li ( \vep, \f{1}{2} \ri )$. \eeL

\bpf  Our computation shows that {\small \[ \f{ \pa
\mscr{M}_{\mrm{B}} (\mu + \vep, \mu) } {\pa \mu}  =  \ln \f{\mu (1 -
\mu - \vep) }{(\mu + \vep) (1 - \mu) }  + \f{ \vep }{\mu (1 - \mu)},
\qqu \f{ \pa^2 \mscr{M}_{\mrm{B}} (\mu + \vep, \mu) } {\pa \mu \pa
\vep} = \f{1}{ \mu (1- \mu)  } -  \f{1}{ (\mu + \vep) (1- \mu -
\vep)  }.
\]} Since $\f{ \pa \mscr{M}_{\mrm{B}} (\mu + \vep, \mu) } {\pa \mu} = 0$ for $\vep = 0$
and $\f{ \pa^2 \mscr{M}_{\mrm{B}} (\mu + \vep, \mu) } {\pa \mu \pa
\vep}
> 0 $ for $\vep < \f{1}{2} - \mu$, it must be true that  $\f{ \pa
\mscr{M}_{\mrm{B}} (\mu + \vep, \mu) } {\pa \mu} > 0$ for $\mu \in
\li (0, \f{1}{2} - \vep \ri )$. This proves the first statement of
the lemma.  Similarly, we can show that {\small $\f{ \pa
\mscr{M}_{\mrm{B}} (\mu - \vep, \mu) } {\pa \mu}  =  \ln \f{\mu (1 -
\mu + \vep) }{(\mu - \vep) (1 - \mu) }  - \f{ \vep }{\mu (1 -
\mu)}$} and {\small $\f{ \pa^2 \mscr{M}_{\mrm{B}} (\mu - \vep, \mu)
} {\pa \mu \pa \vep} = \f{1}{ (\mu - \vep) (1- \mu + \vep)  } -
\f{1}{ \mu (1- \mu) }$}.  Since $\f{ \pa \mscr{M}_{\mrm{B}} (\mu -
\vep, \mu) } {\pa \mu} = 0$ for $\vep = 0$ and $\f{ \pa^2
\mscr{M}_{\mrm{B}} (\mu - \vep, \mu) } {\pa \mu \pa \vep}
> 0 $ for $0 < \vep < \mu < \f{1}{2}$, it must be true that  $\f{
\pa \mscr{M}_{\mrm{B}} (\mu - \vep, \mu) } {\pa \mu} > 0$ for $0 <
\vep < \mu < \f{1}{2}$.  This proves the second statement of the
lemma.

 \epf

 \beL \la{lem8} Let $0 < \vep < \f{1}{2}$. Then,
 $\mscr{M}_{\mrm{B}} (\mu + \vep, \mu) > \mscr{M}_{\mrm{B}} (\mu - \vep, \mu) $
 for $\mu \in \li ( \vep, \f{1}{2} \ri )$. \eeL

\bpf

Straightforward computation shows that
\[
\f{\pa  \mscr{M}_{\mrm{B}} (\mu + \vep, \mu) } {\pa \vep  } =  \ln
\li ( \f{ \mu } { 1 - \mu } \f{1 - \mu - \vep } { \mu + \vep } \ri
), \qqu \f{\pa \mscr{M}_{\mrm{B}} ( \mu - \vep, \mu) } {\pa \vep } =
- \ln \li ( \f{ \mu } { 1 - \mu } \f{1 - \mu + \vep } { \mu - \vep }
\ri ).
\]
Thus, {\small $\f{\pa  \mscr{M}_{\mrm{B}} (\mu + \vep, \mu) } {\pa
\vep } - \f{\pa \mscr{M}_{\mrm{B}} (\mu - \vep, \mu) } {\pa \vep  }
= \ln \f{ \mu^2 } { (1 - \mu)^2 } \f{(1 - \mu)^2 - \vep^2 } { \mu^2
- \vep^2 }
> 0$} if $\vep < \mu < \f{1}{2}$. By virtue of such result and the fact that
$\mscr{M}_{\mrm{B}} (\mu + \vep, \mu)
= \mscr{M}_{\mrm{B}} (\mu - \vep, \mu)$ for $\vep = 0$, we have
$\mscr{M}_{\mrm{B}} (\mu + \vep, \mu)
> \mscr{M}_{\mrm{B}} (\mu - \vep, \mu)$ for $ \vep < \mu < \f{1}{2}$.  This
proves the lemma.

\epf

 \beL \la{lem9} Let $0 < \vep < 1$. Then, $\mscr{M}_{\mrm{B}} (\mu + \vep \mu, \mu)$
 is monotonically decreasing
 with respect to $\mu \in \li (0, \f{1}{1 + \vep} \ri )$.  Similarly,
$\mscr{M}_{\mrm{B}} (\mu - \vep \mu, \mu)$ is monotonically
decreasing with respect to $\mu \in (0, 1)$. \eeL

\bpf

The first statement of the lemma is true because {\small \bee \f{
\pa \mscr{M}_{\mrm{B}} (\mu + \vep \mu, \mu) } {\pa \mu}  = (1 +
\vep) \ln \li [ 1 - \f{ \vep }{(1 + \vep)(1 - \mu)} \ri ] + \f{\vep
}{1 - \mu} < (1 + \vep) \times \li [ - \f{ \vep }{(1 + \vep)(1 -
\mu)} \ri ] + \f{\vep }{1 - \mu} = 0 \eee} for $0 < \mu < \f{1}{1 +
\vep}$. Similarly, the second statement of the lemma is true because
{\small \bee \f{ \pa \mscr{M}_{\mrm{B}} (\mu - \vep \mu, \mu) } {\pa
\mu}  =
 (1 - \vep) \ln \li [ 1 + \f{ \vep }{(1 - \vep)(1 - \mu)} \ri ] - \f{\vep }{1 - \mu}
  <  (1 - \vep) \times
  \li [  \f{ \vep }{(1 - \vep)(1 - \mu)} \ri ] - \f{\vep }{1 - \mu} = 0 \eee} for $0 < \mu < 1$.

\epf

 \beL \la{lem10} Let $0 < \vep < 1$. Then,
 $\mscr{M}_{\mrm{B}} (\mu + \vep \mu, \mu) > \mscr{M}_{\mrm{B}} (\mu - \vep \mu, \mu)$
 for  $\mu \in \li ( 0, \f{1}{2} \ri )$. \eeL

\bpf

It can be shown by tedious computation that \[
 \f{ \pa \mscr{M}_{\mrm{B}} (\mu + \vep \mu, \mu) }{\pa \vep}
 =  \mu \ln \f{ 1 - \mu - \vep \mu } { (1 + \vep) ( 1 - \mu) }, \qqu
 \f{ \pa \mscr{M}_{\mrm{B}} (\mu - \vep \mu, \mu) }{\pa \vep} = - \mu \ln \f{ 1 - \mu + \vep \mu } { (1 -
\vep) ( 1 - \mu) }.
\]
Hence,
\[
\f{ \pa \mscr{M}_{\mrm{B}} (\mu + \vep \mu, \mu) }{\pa \vep} - \f{
\pa \mscr{M}_{\mrm{B}} (\mu - \vep \mu, \mu) }{\pa \vep} = \mu \ln
\li [ \f{ 1 - \mu - \vep \mu } { (1 + \vep) ( 1 - \mu) } \; \f{ 1 -
\mu + \vep \mu } { (1 - \vep) ( 1 - \mu) }  \ri ].
\]
Since {\small $\f{ 1 - \mu - \vep \mu } { (1 + \vep) ( 1 - \mu) } \;
\f{ 1 - \mu + \vep \mu } { (1 - \vep) ( 1 - \mu) } = \f{ ( 1 -
\mu)^2 - \vep^2 \mu^2 } { ( 1 - \mu)^2  - \vep^2 ( 1 - \mu)^2 } >
1$}  for $0 < \mu < \f{1}{2}$, we have {\small $\f{ \pa
\mscr{M}_{\mrm{B}} (\mu + \vep \mu, \mu) }{\pa \vep} - \f{ \pa
\mscr{M}_{\mrm{B}} (\mu - \vep \mu, \mu) }{\pa \vep}
> 0$} for $0 < \mu < \f{1}{2}$. Noting that $\mscr{M}_{\mrm{B}} (\mu
+ \vep \mu, \mu) - \mscr{M}_{\mrm{B}} (\mu - \vep \mu, \mu) = 0$ for
$\vep = 0$, we have $\mscr{M}_{\mrm{B}} (\mu + \vep \mu, \mu) -
\mscr{M}_{\mrm{B}} (\mu - \vep \mu, \mu)
> 0$ for $0 < \mu < \f{1}{2}$. This completes the proof of the lemma.

\epf

\beL  \la{lem11} $\Pr \{ \wh{\bs{p}} \geq (1 + \vep_r) p \} <
\f{\de}{2}$ for any $p \in \li ( p^\star, 1 \ri )$.

 \eeL

 \bpf

To prove the lemma, we shall consider the following three cases:

 Case (i):  $(1 + \vep_r) p > 1$;

 Case (ii):  $\f{\ga}{n} \leq (1 + \vep_r) p \leq 1$;

 Case (iii):  $(1 + \vep_r) p < \f{\ga}{n}$.

\bsk

For Case (i), it is obvious that $\Pr \{ \wh{\bs{p}} \geq (1 +
\vep_r) p \} = 0 < \f{\de}{2}$.

For Case (ii), applying Theorem 1 with $z = (1 + \vep_r) p \geq
\f{\ga}{n}$, we have \bel \Pr \{ \wh{\bs{p}} \geq (1 + \vep_r) p \}
& = & \Pr \li \{ \sum_{i = 1}^{\lf \ga \sh z \rf} X_i \geq \ga \ri
\} \nonumber\\
& \leq & \exp \li (  \lf \ga \sh z \rf \; \mscr{M}_{\mrm{B}} \li (
\f{\ga}{\lf \ga \sh z \rf}, p
 \ri ) \ri ) \la{uselem5} \\
 & = & \exp \li (  \ga \; \mscr{M}_{\mrm{I}} \li (
\f{\ga}{\lf \ga \sh z \rf}, p
 \ri ) \ri ) \nonumber\\
 & \leq & \exp \li (  \ga \; \mscr{M}_{\mrm{I}} \li ( z, p
 \ri ) \ri ) \la{monde}\\
 & < & \exp \li ( \ga \mscr{M}_{\mrm{I}} ( p^\star + \vep_r p^\star, p^\star  ) \ri ) \la{e121}\\
& < & \f{\de}{2} \la{e122} \eel where (\ref{uselem5}) follows from
Lemma \ref{lem1},  (\ref{monde}) is due to the fact that
$\mscr{M}_{\mrm{I}} \li ( z, p \ri )$ is monotonically decreasing
with respect to $z \in (p, 1)$, (\ref{e121}) follows from Lemma
\ref{lem3}, and (\ref{e122}) follows from  the assumption about
$\ga$.

For Case (iii), applying Theorem 1 with $z = (1 + \vep_r) p <
\f{\ga}{n}$, we have \bel \Pr \{ \wh{\bs{p}} \geq (1 + \vep_r) p \}
& = &  \Pr \li \{ \sum_{i = 1}^{n} X_i \geq n z
 \ri \} \nonumber\\
& \leq & \exp \li ( n  \mscr{M}_{\mrm{B}}( p + \vep_r p, p  ) \ri ) \la{e12a}\\
& < & \exp \li ( n  \mscr{M}_{\mrm{B}}( p^\star + \vep_r p^\star, p^\star  ) \ri ) \la{e12b}\\
& = & \exp \li ( n  \mscr{M}_{\mrm{B}}( p^\star + \vep_a, p^\star  ) \ri ) \nonumber\\
& < & \f{\de}{2} \la{e12c} \eel where (\ref{e12a}) follows from
Lemma \ref{lem1}, (\ref{e12b}) follows from the first statement of
Lemma \ref{lem9}, and (\ref{e12c}) follows from the assumption about
$n$.

  In summary, we have shown $\Pr \{ \wh{\bs{p}} \geq (1
+ \vep_r) p \} < \f{\de}{2}$ for all cases. This completes the proof
of the lemma.

\epf

\beL  \la{lem12}  $\Pr \{ \wh{\bs{p}} \leq (1 - \vep_r) p \} <
\f{\de}{2}$ for any $p \in \li ( p^\star, 1 \ri )$. \eeL

\bpf To prove the lemma, we shall consider the following two cases:

Case (i): $(1 - \vep_r) p \geq \f{\ga}{n}$;

Case (ii): $(1 - \vep_r) p < \f{\ga}{n}$.

\bsk

For Case (i),  applying Theorem 1 with $z = (1 - \vep_r) p \geq
\f{\ga}{n}$, we have {\small \bel \Pr \{ \wh{\bs{p}} \leq (1 -
\vep_r) p \} & = & \Pr \li \{ \sum_{i = 1}^{\lc \ga \sh z \rc - 1}
X_i < \ga \ri \} =  \Pr \li \{ \sum_{i = 1}^{\lc \ga \sh z \rc - 1}
X_i \leq  \ga - 1 \ri \} \leq  \Pr \li \{ \sum_{i = 1}^{\lc \ga \sh
z \rc} X_i \leq \ga \ri \} \nonumber\\
& \leq & \exp \li (  \lc \ga \sh z \rc \; \mscr{M}_{\mrm{B}} \li (
\f{\ga}{\lc \ga \sh z \rc}, p \ri ) \ri ) \la{useH}\\
 & = & \exp \li (  \ga \; \mscr{M}_{\mrm{I}} \li (
\f{\ga}{\lc \ga \sh z \rc}, p
 \ri ) \ri ) \nonumber\\
 & \leq & \exp \li (  \ga \; \mscr{M}_{\mrm{I}} \li ( z, p
 \ri ) \ri ) \la{dez}\\
 & < & \exp \li ( \ga \mscr{M}_{\mrm{I}} ( p^\star - \vep_r p^\star, p^\star  ) \ri ) \la{e1218}\\
& < & \exp \li ( \ga \mscr{M}_{\mrm{I}} ( p^\star + \vep_r p^\star,
p^\star  ) \ri ) \la{newlem} \\
&  < & \f{\de}{2}, \nonumber \eel} where (\ref{useH}) follows from
Lemma \ref{lem1},  (\ref{dez}) is due to the fact that
$\mscr{M}_{\mrm{I}} \li ( z, p \ri )$ is monotonically increasing
with respect to $z \in (0, p)$,  (\ref{e1218}) follows from Lemma
\ref{lem3}, and (\ref{newlem}) follows from Lemma \ref{lem4}.

For Case (ii), applying Theorem 1 with $z = (1 - \vep_r) p <
\f{\ga}{n}$, we have \bel\Pr \{ \wh{\bs{p}} \leq (1 - \vep_r) p \} &
\leq & \exp \li ( n \mscr{M}_{\mrm{B}} ( p - \vep_r p, p  ) \ri ) \la{c12a}\\
& < & \exp \li ( n  \mscr{M}_{\mrm{B}} ( p^\star - \vep_r p^\star, p^\star  ) \ri ) \la{c12b}\\
& < & \exp \li ( n  \mscr{M}_{\mrm{B}} ( p^\star + \vep_r p^\star, p^\star  ) \ri ) \la{c12c}\\
& = &  \exp \li ( n  \mscr{M}_{\mrm{B}} ( p^\star + \vep_a, p^\star
) \ri ) < \f{\de}{2} \nonumber \eel where (\ref{c12a}) follows from
Lemma \ref{lem1}, (\ref{c12b}) follows from the second statement of
Lemma \ref{lem9}, and (\ref{c12c}) follows from Lemma \ref{lem10}.

In summary, we have shown $\Pr \{ \wh{\bs{p}} \leq (1 - \vep_r) p \}
< \f{\de}{2}$ for both cases.  The lemma is thus proved.

\epf

\beL  \la{lem13}  $\Pr \{ \wh{\bs{p}} \geq p  + \vep_a \} <
\f{\de}{2}$ for any $p \in \li (0, p^\star \ri ]$. \eeL

\bpf  To prove the lemma, we shall consider the following two cases:

Case (i): $p + \vep_a \geq \f{\ga}{n}$;

Case (ii): $p + \vep_a < \f{\ga}{n}$.

\bsk

For Case (i), applying Theorem 1 with $z = p + \vep_a \geq
\f{\ga}{n}$, we have {\small \bel \Pr \{ \wh{\bs{p}} \geq p + \vep_a
\} & = & \Pr \li \{ \sum_{i = 1}^{\lf \ga \sh z \rf} X_i \geq \ga
\ri \} \leq \exp \li ( \lf \ga \sh z \rf \; \mscr{M}_{\mrm{B}} \li (
\f{\ga}{\lf \ga \sh z \rf}, p
 \ri ) \ri ) \nonumber\\
 & = & \exp \li (  \ga \; \mscr{M}_{\mrm{I}} \li (
\f{\ga}{\lf \ga \sh z \rf}, p
 \ri ) \ri ) \leq \exp \li (  \ga \; \mscr{M}_{\mrm{I}} \li ( z, p
 \ri ) \ri ) \nonumber\\
 & < & \exp \li ( \ga \mscr{M}_{\mrm{I}} ( p^\star + \vep_r p^\star, p^\star  ) \ri ) \la{e1219}\\
& < & \f{\de}{2}, \nonumber \eel} where (\ref{e1219}) follows from
Lemma \ref{lem5}.

For Case (ii), applying Theorem 1 with $z = p + \vep_a <
\f{\ga}{n}$, we have \bel \Pr \{ \wh{\bs{p}} \geq p + \vep_a \} & =
& \Pr \li \{ \sum_{i =
1}^{n} X_i \geq n z  \ri \} \nonumber\\
& \leq & \exp \li ( n \mscr{M}_{\mrm{B}} ( p + \vep_a, p  ) \ri ) \la{c13c}\\
& \leq & \exp \li ( n \mscr{M}_{\mrm{B}} ( p^\star + \vep_a, p^\star  ) \ri ) \la{c13d}\\
& < & \f{\de}{2} \nonumber \eel where (\ref{c13c}) follows from
Lemma \ref{lem1}, (\ref{c13d}) follows from the first statement of
Lemma \ref{lem7}.

In summary, we have shown $\Pr \{ \wh{\bs{p}} \geq p  + \vep_a \} <
\f{\de}{2}$ for both cases.  The lemma is thus proved.

\epf

\beL

\la{lem14}

$\Pr \{ \wh{\bs{p}} \leq p - \vep_a \} < \f{\de}{2}$ for any $p \in
\li (0,  p^\star \ri ]$.

\eeL

\bpf

To prove the lemma, we shall consider the following three cases:

Case (i): $p < \vep_a$;

Case (ii): $p - \vep_a \geq \f{\ga}{n}$;

Case (iii): $0 \leq p - \vep_a < \f{\ga}{n}$.

\bsk

For Case (i), it is obvious that $\Pr \{ \wh{\bs{p}} \leq p - \vep_a
\} = 0 < \f{\de}{2}$.

For Case (ii), applying Theorem 1 with $z = p - \vep_a \geq
\f{\ga}{n}$, we have \bel \Pr \{ \wh{\bs{p}} \leq p - \vep_a \} & =
& \Pr \li \{ \sum_{i = 1}^{\lc \ga \sh z \rc - 1} X_i < \ga \ri \}
\leq  \Pr \li \{ \sum_{i = 1}^{\lc \ga \sh z \rc} X_i \leq \ga \ri
\} \leq \exp \li (  \lc \ga \sh z \rc \; \mscr{M}_{\mrm{B}} \li (
\f{\ga}{\lc \ga \sh z \rc}, p
 \ri ) \ri ) \nonumber \\
 &  =  & \exp \li (  \ga \; \mscr{M}_{\mrm{I}} \li (
\f{\ga}{\lc \ga \sh z \rc}, p
 \ri ) \ri ) \leq \exp \li (  \ga \; \mscr{M}_{\mrm{I}} \li ( z, p
 \ri ) \ri )  =  \exp \li ( \ga \mscr{M}_{\mrm{I}} ( p - \vep_a, p  ) \ri ) \nonumber \\
& \leq & \exp \li ( \ga \mscr{M}_{\mrm{I}} ( p^\star - \vep_a, p^\star  ) \ri ) \la{c14b}\\
& \leq & \exp \li ( \ga \mscr{M}_{\mrm{I}} ( p^\star + \vep_a, p^\star  ) \ri ) \la{c14c}\\
&  < & \f{\de}{2} \nonumber \eel where  (\ref{c14b}) follows from
the second statement of Lemma \ref{lem5}, and (\ref{c14c}) follows
from Lemma \ref{lem6}.

For Case (iii), applying Theorem 1 with $z = p - \vep_a <
\f{\ga}{n}$, we have \bel \Pr \{ \wh{\bs{p}} \leq p - \vep_a \} & =
& \Pr \li \{ \sum_{i = 1}^n X_i \leq n z
\ri \} \nonumber\\
& \leq & \exp \li ( n \mscr{M}_{\mrm{B}} ( p - \vep_a, p  ) \ri ) \la{c14d}\\
& \leq & \exp \li ( n  \mscr{M}_{\mrm{B}} ( p^\star - \vep_a, p^\star  ) \ri ) \la{c14e}\\
& \leq & \exp \li ( n  \mscr{M}_{\mrm{B}} ( p^\star + \vep_a, p^\star  ) \ri ) \la{c14f}\\
& < & \f{\de}{2} \nonumber \eel where (\ref{c14d}) follows from
Lemma \ref{lem1}, (\ref{c14e}) follows from the second statement of
Lemma \ref{lem7}, and (\ref{c14f}) follows from Lemma \ref{lem8}.

In summary, we have shown $\Pr \{ \wh{\bs{p}} \leq p - \vep_a \} < \f{\de}{2}$ for all cases.
The lemma is thus proved.

 \epf

Now we are in a position to prove Theorem 3.  To show $\Pr \li \{
\li | \wh{\bs{p}} - p \ri | < \vep_a  \; \tx{or} \; \li |
\wh{\bs{p}} - p \ri | < \vep_r p \ri \} > 1 - \de$, it suffices to
show $\Pr \li \{  \li | \wh{\bs{p}} - p \ri | \geq \vep_a, \;  \li |
\wh{\bs{p}} - p \ri | \geq \vep_r p \ri \} < \de$ for $0 < p < 1$.

For $p \in \li ( p^\star, 1 \ri )$, we have \bel \Pr \li \{  \li |
\wh{\bs{p}} - p \ri | \geq \vep_a, \;  \li | \wh{\bs{p}} - p \ri |
\geq \vep_r p
\ri \} & = & \Pr \li \{  \li | \wh{\bs{p}} - p \ri | \geq \vep_r p \ri \}
 =  \Pr \{ \wh{\bs{p}} \geq (1 + \vep_r) p \} +
 \Pr \{ \wh{\bs{p}} \leq (1 - \vep_r) p \} \nonumber\\
& < & \f{\de}{2} + \f{\de}{2} \la{ccc}\\
& = & \de \nonumber \eel where (\ref{ccc}) follows from Lemmas
\ref{lem11} and \ref{lem12}.  Similarly, for $p \in (0,  p^\star]$,
we have \bel \Pr \li \{  \li | \wh{\bs{p}} - p \ri | \geq \vep_a, \;
\li | \wh{\bs{p}} - p \ri | \geq \vep_r p \ri \} & = & \Pr \li \{
\li | \wh{\bs{p}} - p \ri | \geq \vep_a \ri \}
 =  \Pr \{ \wh{\bs{p}} \geq p + \vep_a \} + \Pr \{ \wh{\bs{p}} \leq p - \vep_a \} \nonumber\\
& < & \f{\de}{2} + \f{\de}{2} \la{cccc}\\
& = & \de \nonumber \eel where (\ref{cccc}) follows from Lemmas
\ref{lem13} and \ref{lem14}.  This completes the proof of Theorem 3.

\sect{Proof of Theorem 4}

\beL \la{constant}  Let $I$ denote the support of $\wh{\bs{p}}$.
Suppose the intersection between open interval $(p^\prime, \;
p^{\prime \prime})$ and  set $I_\mscr{L}$ is empty.  Then, $\{ \vse
\in I : p \leq \mscr{L} (\vse) \}$ is fixed with respect to $p \in
(p^\prime, \; p^{\prime \prime})$. \eeL

\bpf

Let $p^*$ and $p^\diamond$ be two distinct real numbers included in
interval $(p^\prime, \; p^{\prime \prime})$.  To show the lemma, it
suffices to show that $\{ \vse \in I : p^* \leq \mscr{L} (\vse) \} =
\{ \vse \in I : p^\diamond \leq \mscr{L} (\vse) \}$. First, we shall
show that $\{ \vse \in I : p^* \leq \mscr{L} (\vse) \} \subseteq \{
\vse \in I : p^\diamond \leq \mscr{L} (\vse) \}$. To this end, we
let $\vpi \in \{ \vse \in I : p^* \leq \mscr{L} (\vse) \}$ and
proceed to show $\vpi \in \{ \vse \in I : p^\diamond \leq \mscr{L}
(\vse) \}$. Since $\vpi \in I$ and $p^* \leq \mscr{L} (\vpi)$, it
must be true that $\vpi \in I$ and $p^\diamond \leq \mscr{L}
(\vpi)$.  If this is not the case, then we have $p^{\prime \prime} >
p^\diamond
> \mscr{L} (\vpi) \geq p^* > p^\prime$. Consequently, $\mscr{L} (\vpi)$
is included by both the interval
$(p^\prime, \; p^{\prime \prime})$ and the set $I_\mscr{L}$. This
contradicts the assumption of the lemma. Hence, we have shown $\vpi
\in \{ \vse \in I : p^\diamond \leq \mscr{L} (\vse) \}$ and
accordingly $\{ \vse \in I : p^* \leq \mscr{L} (\vse) \} \subseteq
\{ \vse \in I : p^\diamond \leq \mscr{L} (\vse) \}$. Second, by a
similar argument, we can show $\{ \vse \in I : p^\diamond \leq
\mscr{L} (\vse) \} \subseteq \{ \vse \in I : p^* \leq \mscr{L}
(\vse) \}$. It follows that $\{ \vse \in I : p^* \leq \mscr{L}
(\vse) \} = \{ \vse \in I : p^\diamond \leq \mscr{L} (\vse) \}$.
Finally, the proof of the lemma is completed by noting that the
above argument holds for arbitrary $p^*$ and $p^\diamond$ included
in the open interval $(p^\prime, \; p^{\prime \prime})$.

\epf

By virtue of Theorem 1, we can show the following lemma.

 \beL
 \la{incom383} $\Pr \{ \wh{\bs{p}} \leq  z \mid p \}$ is monotonically
decreasing with respect to $p$. Similarly, $\Pr \{ \wh{\bs{p}} \geq
z \mid p \}$ is monotonically increasing with respect to $p$.

\eeL

\beL \la{sim8} Let $p^\prime < p^{\prime \prime}$ be two consecutive
distinct elements of $I_\mscr{L} \cap [a, b] \cup \{ a, b\}$.  Then,
\bee & & \lim_{\ep \downarrow 0} \Pr \{ p^\prime + \ep \leq
\mscr{L}( \wh{\bs{p}} ) \mid p^\prime + \ep \}  = \Pr \{
p^\prime < \mscr{L}( \wh{\bs{p}} ) \mid p^\prime \},\\
&  & \lim_{\ep \downarrow 0} \Pr \{ p^{\prime \prime} - \ep \leq
\mscr{L}( \wh{\bs{p}} ) \mid p^{\prime \prime} - \ep \}  =  \Pr \{
p^{\prime \prime} \leq \mscr{L}( \wh{\bs{p}} ) \mid p^{\prime
\prime} \}. \eee Moreover,  $\Pr \{ p \leq  \mscr{L}( \wh{\bs{p}} )
\mid p \}$ is monotone with respect to $p \in (p^\prime, p^{\prime
\prime})$. \eeL

\bpf

First, we shall show that $\lim_{\ep \downarrow 0} \Pr \{ p^\prime +
\ep \leq \mscr{L}( \wh{\bs{p}} ) \mid p^\prime + \ep \}  = \Pr \{
p^\prime < \mscr{L}( \wh{\bs{p}} ) \mid p^\prime \}$.   Let
$m^+(\ep)$ be the number of elements of $\{ \vse \in I : p^\prime <
\mscr{L}( \vse ) < p^\prime + \ep
 \}$, where $I$ denotes the support of $\wh{\bs{p}}$ as in Lemma \ref{constant}.   We claim
that $\lim_{\ep \downarrow 0} m^+(\ep) = 0$. It suffices to consider
two cases as follows.

In the case of $\{ \vse \in I : p^\prime < \mscr{L}( \vse ) \} =
\emptyset$, we have $m^+(\ep) = 0$ for any $\ep >
 0$.  In the case of $\{ \vse \in I : p^\prime < \mscr{L}( \vse ) \} \neq \emptyset$,
 we have $m^+(\ep) = 0$ for
 $0 < \ep \leq \ep^*$, where $\ep^* = \min \li \{  \mscr{L}( \vse ) - p^\prime : p^\prime <
\mscr{L}( \vse ), \; \vse \in I \ri \}$ is positive because of the
assumption that $I_\mscr{L}$ has no closure points in $[a, b]$.
Hence, in both cases, $\lim_{\ep \downarrow 0} m^+(\ep) = 0$. This
establishes the claim.

Noting that $\Pr \{ p^\prime  < \mscr{L}( \wh{\bs{p}} ) < p^\prime +
\ep \mid p^\prime + \ep \} \leq m^+(\ep)$ as a consequence of $\Pr
\{ \wh{\bs{p}} = \vse \mid p^\prime + \ep \} \leq 1$ for any $\vse
\in I$, we have that $\limsup_{\ep \downarrow 0} \Pr \{ p^\prime  <
\mscr{L}( \wh{\bs{p}} ) < p^\prime + \ep \mid p^\prime + \ep \} \leq
\lim_{\ep \downarrow 0} m^+(\ep) = 0$, which implies that $\lim_{\ep
\downarrow 0} \Pr \{ p^\prime  < \mscr{L}( \wh{\bs{p}} ) < p^\prime
+ \ep \mid p^\prime + \ep \} = 0$.

Since $\{ p^\prime + \ep \leq \mscr{L}( \wh{\bs{p}} ) \} \cap \{
p^\prime  < \mscr{L}( \wh{\bs{p}} ) < p^\prime + \ep \} = \emptyset$
and $\{ p^\prime  < \mscr{L}( \wh{\bs{p}} ) \}
 = \{ p^\prime + \ep \leq \mscr{L}( \wh{\bs{p}} ) \} \cup
 \{ p^\prime  < \mscr{L}( \wh{\bs{p}} ) <
p^\prime + \ep \}$,  we have $\Pr \{ p^\prime < \mscr{L}(
\wh{\bs{p}} ) \mid p^\prime + \ep \} = \Pr \{ p^\prime + \ep \leq
\mscr{L}( \wh{\bs{p}} ) \mid p^\prime + \ep \} + \Pr \{ p^\prime <
\mscr{L}( \wh{\bs{p}} ) < p^\prime + \ep \mid p^\prime + \ep \}$.
Observing that $\Pr \{ p^\prime < \mscr{L}( \wh{\bs{p}} ) \mid
p^\prime + \ep \}$ is continuous with respect to $\ep \in (0, 1 -
p^\prime )$, we have $\lim_{\ep \downarrow 0} \Pr \{ p^\prime  <
\mscr{L}( \wh{\bs{p}} ) \mid p^\prime + \ep \} = \Pr \{ p^\prime  <
\mscr{L}( \wh{\bs{p}} ) \mid p^\prime \}$. It follows that {\small
\bee \lim_{\ep \downarrow 0} \Pr \{ p^\prime + \ep \leq \mscr{L}(
\wh{\bs{p}} ) \mid p^\prime + \ep \} & = & \lim_{\ep \downarrow 0}
\Pr \{ p^\prime < \mscr{L}( \wh{\bs{p}} ) \mid p^\prime + \ep \} -
\lim_{\ep \downarrow 0} \Pr \{ p^\prime < \mscr{L}( \wh{\bs{p}}
) < p^\prime + \ep \mid p^\prime + \ep \}\\
& = & \lim_{\ep \downarrow 0} \Pr \{ p^\prime < \mscr{L}(
\wh{\bs{p}} ) \mid p^\prime + \ep \} =  \Pr \{ p^\prime  < \mscr{L}(
\wh{\bs{p}} ) \mid p^\prime \}. \eee}

Next, we shall show that $\lim_{\ep \downarrow 0} \Pr \{ p^{\prime
\prime} - \ep \leq \mscr{L}( \wh{\bs{p}} ) \mid p^{\prime \prime} -
\ep \} = \Pr \{ p^{\prime \prime} \leq \mscr{L}( \wh{\bs{p}} ) \mid
p^{\prime \prime} \}$.   Let $m^{-} (\ep) $ be the number of
elements of $\{ \vse \in I : p^{\prime \prime} - \ep \leq \mscr{L}(
\vse ) < p^{\prime \prime} \}$.  Then, we can show $\lim_{\ep
\downarrow 0} m^{-} (\ep) = 0$ by considering two cases as follows.

In the case of $\{ \vse \in I : \mscr{L}( \vse ) < p^{\prime \prime}
\} = \emptyset$, we have $m^{-} (\ep) = 0$ for any $\ep > 0$. In the
case of $\{ \vse \in I : \mscr{L}( \vse ) < p^{\prime \prime} \}
\neq \emptyset$, we have $m^{-} (\ep) = 0$ for $0 < \ep <
\ep^\star$, where $\ep^\star =  \min \{ p^{\prime \prime} -
\mscr{L}( \vse ) : \vse \in I, \; \mscr{L}( \vse ) < p^{\prime
\prime} \}$ is positive because of the assumption that $I_\mscr{U}$
has no closure points in $[a, b]$.  Hence, in both cases, $\lim_{\ep
\downarrow 0} m^{-} (\ep) = 0$. It follows that $\limsup_{\ep
\downarrow 0} \Pr \{ p^{\prime \prime} - \ep \leq \mscr{L}(
\wh{\bs{p}} ) < p^{\prime \prime} \mid p^{\prime \prime} - \ep \}
\leq \lim_{\ep \downarrow 0} m^{-} (\ep) = 0$ and consequently
$\lim_{\ep \downarrow 0} \Pr \{ p^{\prime \prime} - \ep \leq
\mscr{L}( \wh{\bs{p}} ) < p^{\prime \prime} \mid p^{\prime \prime} -
\ep \} = 0$.

Since $\{ p^{\prime \prime} - \ep \leq \mscr{L}( \wh{\bs{p}} ) \} =
\{ p^{\prime \prime}  \leq \mscr{L}( \wh{\bs{p}} ) \} \cup \{
p^{\prime \prime} - \ep \leq \mscr{L}( \wh{\bs{p}} ) < p^{\prime
\prime} \}$ and $\{ p^{\prime \prime}  \leq \mscr{L}( \wh{\bs{p}} )
\} \cap \{ p^{\prime \prime} - \ep \leq \mscr{L}( \wh{\bs{p}} ) <
p^{\prime \prime} \} = \emptyset$, we have $\Pr \{ p^{\prime \prime}
- \ep \leq \mscr{L}( \wh{\bs{p}} ) \mid p^{\prime \prime} - \ep \} =
\Pr \{ p^{\prime \prime}  \leq \mscr{L}( \wh{\bs{p}} ) \mid
p^{\prime \prime} - \ep \} + \Pr \{ p^{\prime \prime} - \ep \leq
\mscr{L}( \wh{\bs{p}} ) < p^{\prime \prime} \mid p^{\prime \prime} -
\ep \}$.

Observing that $\Pr \{ p^{\prime \prime} \leq \mscr{L}( \wh{\bs{p}}
) \mid p^{\prime \prime} - \ep \}$ is continuous with respect to
$\ep \in (0, p^{\prime \prime})$,  we have $\lim_{\ep \downarrow 0}
\Pr \{ p^{\prime \prime} \leq \mscr{L}( \wh{\bs{p}} ) \mid p^{\prime
\prime} - \ep \} = \Pr \{ p^{\prime \prime}  \leq \mscr{L}(
\wh{\bs{p}} ) \mid p^{\prime \prime} \}$. It follows that $\lim_{\ep
\downarrow 0} \Pr \{ p^{\prime \prime} - \ep \leq \mscr{L}(
\wh{\bs{p}} ) \mid p^{\prime \prime} - \ep \} = \lim_{\ep \downarrow
0} \{ p^{\prime \prime} \leq \mscr{L}( \wh{\bs{p}} ) \mid p^{\prime
\prime} \}$.

Now we turn to show that $\Pr \{ p \leq  \mscr{L}( \wh{\bs{p}} )
\mid p \}$ is monotone with respect to $p \in (p^\prime, p^{\prime
\prime})$.   Without loss of generality, we assume that
$\mscr{L}(.)$ is monotonically increasing.  Since $p^\prime <
p^{\prime \prime}$ are two consecutive distinct elements of
$I_\mscr{L} \cap [a, b] \cup \{ a, b\}$, we have that the
intersection between open interval $(p^\prime, p^{\prime \prime})$
and  set $I_\mscr{L}$ is empty.  As a result of Lemma
\ref{constant}, we can write $\Pr \{ p \leq \mscr{L}( \wh{\bs{p}} )
\mid p \} = \Pr \{ \wh{\bs{p}} \geq \vse \mid p \}$, where $\vse \in
[0, 1]$ is a constant independent of $p \in (p^\prime, p^{\prime
\prime})$. By Lemma \ref{incom383}, we have that $\Pr \{ \wh{\bs{p}}
\geq \vse \mid p \}$ is monotonically increasing with respect to $p
\in (p^\prime, p^{\prime \prime})$. This proves the monotonicity of
$\Pr \{ p \leq \mscr{L}( \wh{\bs{p}} ) \mid p \}$ with respect to $p
\in (p^\prime, p^{\prime \prime})$. The proof of the lemma is thus
completed.

\epf

By a similar method as that of Lemma \ref{sim8}, we can show the
following lemma.

\beL Let $p^\prime < p^{\prime \prime}$ be two consecutive distinct
elements of $I_\mscr{U} \cap [a, b] \cup \{ a, b\}$.  Then, \bee  &
& \lim_{\ep \downarrow 0} \Pr \{ p^\prime + \ep \geq \mscr{U}(
\wh{\bs{p}} ) \mid p^\prime + \ep \} = \Pr \{ p^\prime \geq
\mscr{U}( \wh{\bs{p}}
) \mid p^\prime \},\\
&  & \lim_{\ep \downarrow 0} \Pr \{ p^{\prime \prime} - \ep \geq
\mscr{U}( \wh{\bs{p}} ) \mid p^{\prime \prime} - \ep \} = \Pr \{
p^{\prime \prime} > \mscr{U}( \wh{\bs{p}} ) \mid p^{\prime \prime}
\}. \eee Moreover,  $\Pr \{ p \geq  \mscr{U}( \wh{\bs{p}} ) \mid p
\}$ is monotone with respect to $p \in (p^\prime, p^{\prime
\prime})$.
 \eeL

Now we are in a position to prove Theorem \ref{THM1}.  Let $C(p) =
\Pr \{ p \leq \mscr{L}( \wh{\bs{p}} ) \mid p \}$.  By Lemma
\ref{sim8}, $C(p)$ is a monotone function of $p \in (p^\prime,
p^{\prime \prime})$, which implies that $C(p) \leq \max \{
C(p^\prime + \ep), \; C(p^{\prime \prime} - \ep) \}$ for any $p \in
(p^\prime, p^{\prime \prime})$ and any positive $\ep$ less than
$\min \{ p - p^\prime, \; p^{\prime \prime} - p \}$. Consequently,
\[
C(p) \leq \lim_{\ep \downarrow 0} \max \{ C(p^\prime + \ep), \;
C(p^{\prime \prime} - \ep) \} = \max \{ \lim_{\ep \downarrow 0}
C(p^\prime + \ep), \; \lim_{\ep \downarrow 0} C(p^{\prime \prime} -
\ep)
 \} \leq \max \{ C(p^\prime), \; C(p^{\prime \prime}) \}
\]
for any $p \in (p^\prime, p^{\prime \prime})$.  Since the argument
holds for arbitrary consecutive distinct elements of {\small $\{
\mscr{L} ( \wh{p} ) \in (a, b) \mid \wh{p} \in I \} \cup \{ a,
b\}$}, we have established the statement regarding the maximum of
$\Pr \{ p \leq \mscr{L}( \wh{\bs{p}} ) \mid p \}$ with respect to $p
\in (a, b)$. By a similar method, we can prove the statement
regarding the maximum of $\Pr \{ p \geq \mscr{U}( \wh{\bs{p}} ) \mid
p \}$ with respect to $p \in (a, b)$. This concludes the proof of
Theorem \ref{THM1}.

\sect{Proof of Theorem 6}

The theorem can be established by showing the following lemmas.

\beL $\Pr \{ p \geq \ovl{\bs{p}} \} \leq \f{\de}{2}$. \eeL

\bpf

By Theorem 1, {\small \[ \Pr \{ \wh{\bs{p}} \leq z \} = \bec \Pr \{
\sum_{i = 1}^{\lc \ga \sh z \rc - 1} X_i <
\ga \} & \tx{for} \; \ga \leq n z,\\
\Pr \li \{ \sum_{i = 1}^n X_i \leq n z \ri \} & \tx{for} \; \ga
> n z. \eec
\]}
Since $X_i$ must be either $0$ or $1$ and $\ga$ is an integer, we
have {\small $\Pr \{ \sum_{i = 1}^{\lc \ga \sh z \rc - 1} X_i < \ga
\} \leq \Pr \{ \sum_{i = 1}^{\lc \ga \sh z \rc} X_i \leq \ga \} $}.
Hence, {\small
\[ \Pr \{ \wh{\bs{p}} \leq z \} \leq \bec \Pr \{ \sum_{i = 1}^{\lc
\ga \sh z \rc} X_i \leq
\ga \} & \tx{for} \; \ga \leq n z,\\
\Pr \li \{ \sum_{i = 1}^n X_i \leq n z \ri \} & \tx{for} \; \ga
> n z. \eec
\]}
Since $X_1, X_2, \cd$ are i.i.d. Bernoulli random variables, we have
$\Pr \{  \wh{\bs{p}} \leq z \} \leq \mscr{G}( z, p)$, where
\[
\mscr{G}( z, p) = \bec \sum_{i = 0}^{\ga} \bi{\lc \ga \sh z \rc}{i}
p^i (1 - p)^{\lc \ga
\sh z \rc - i}  &  \tx{for} \; \f{\ga}{n} \leq z \leq 1,\\
\sum_{i = 0}^{\lf n z \rf}  \bi{n}{i} p^i (1 - p)^{n - i}  &
\tx{for} \; 0 \leq z < \f{\ga}{n}. \eec
\]
Let $z^* \in [0, 1]$ be the largest number such that $\Pr \{
\wh{\bs{p}} < z^* \} \leq \f{\de}{2}$.  Since $\wh{\bs{p}}$ is a
discrete random variable bounded in $[0, 1]$, it must be true that
$\Pr \{ \wh{\bs{p}} \leq z^* \} > \f{\de}{2}$. Observing that
$\mscr{G} (z, p)$ is monotonically decreasing with respect to $p \in
(0, 1)$, we have {\small \[ \{ p \geq \ovl{\bs{p}} \} = \{ p \geq
\ovl{\bs{p}}, \; \mbf{k} < \mbf{n} \} \subseteq  \li \{ \mscr{G}
(\wh{\bs{p}}, p) \leq \mscr{G} (\wh{\bs{p}}, \ovl{\bs{p}}) =
\f{\de}{2} \ri \} \subseteq \li \{ \mscr{G} (\wh{\bs{p}}, p) \leq
\f{\de}{2} \ri \}.
\]}
Noting that $\f{\de}{2} < \Pr \{ \wh{\bs{p}} \leq z^* \} \leq
\mscr{G} (z^*, p)$ and that $\mscr{G} (z, p)$ is non-decreasing with
respect to $z \in (0, 1)$, we have $\{ p \geq \ovl{\bs{p}} \}
\subseteq \{ \mscr{G} (\wh{\bs{p}}, p) \leq \f{\de}{2} \} \subseteq
\{ \mscr{G} (\wh{\bs{p}}, p)  < \mscr{G} (z^*, p) \} \subseteq \{
\wh{\bs{p}} < z^* \}$.  It follows that $\Pr \{ p \geq \ovl{\bs{p}}
\} \leq \Pr \{ \wh{\bs{p}} < z^* \} \leq \f{\de}{2}$.

\epf

\beL $\Pr \{ p \leq \udl{\bs{p}} \} \leq \f{\de}{2}$. \eeL

\bpf

 By Theorem 1, {\small
\[
\Pr \{  \wh{\bs{p}} \geq z \} = \bec \Pr  \{ \sum_{i = 1}^{ \lf \ga
\sh z \rf} X_i
\geq \ga  \} & \tx{for} \; \ga \leq n z,\\
\Pr \{ \sum_{i = 1}^{n} X_i \geq n z  \} & \tx{for} \; \ga
> n z. \eec
\]}
Since $X_1, X_2, \cd$ are i.i.d. Bernoulli random variables, we have
$\Pr \{  \wh{\bs{p}} \geq z \} = \mscr{H}( z, p)$ where
\[
\mscr{H}( z, p) = \bec \sum_{i = \ga}^{\lf \ga \sh z \rf} \bi{\lf
\ga \sh z \rf}{i} p^i (1 - p)^{\lf \ga
\sh z \rf - i}   &  \tx{for} \; \f{\ga}{n} \leq z \leq 1,\\
\sum_{i = \lc n z \rc}^n \bi{n}{i} p^i (1 - p)^{n - i}  &  \tx{for}
\; 0 \leq z < \f{\ga}{n}. \eec
\]
Let $z^* \in [0, 1]$ be the smallest number such that $\Pr \{
\wh{\bs{p}} > z^* \} \leq \f{\de}{2}$.  Since $\wh{\bs{p}}$ is a
discrete random variable bounded in $[0, 1]$, it must be true that
$\Pr \{ \wh{\bs{p}} \geq z^* \} > \f{\de}{2}$.  Observing that
$\mscr{H} (z, p)$ is monotonically increasing with respect to $p \in
(0, 1)$, we have {\small \[ \{ p \leq \udl{\bs{p}} \} = \{ p \leq
\udl{\bs{p}}, \; \mbf{k} > 0 \} \subseteq  \li \{ \mscr{H}
(\wh{\bs{p}}, p) \leq \mscr{H} (\wh{\bs{p}}, \udl{\bs{p}}) =
\f{\de}{2} \ri \} \subseteq \li \{ \mscr{H} (\wh{\bs{p}}, p) \leq
\f{\de}{2} \ri \}.
\]}
Noting that $\f{\de}{2} < \Pr \{ \wh{\bs{p}} \geq z^* \} = \mscr{H}
(z^*, p)$ and that $\mscr{H} (z, p)$ is non-increasing with respect
to $z \in (0, 1)$, we have $\{ p \leq \udl{\bs{p}} \} \subseteq  \{
\mscr{H} (\wh{\bs{p}}, p) \leq \f{\de}{2} \} \subseteq \{ \mscr{H}
(\wh{\bs{p}}, p)  < \mscr{H} (z^*, p) \} \subseteq \{ \wh{\bs{p}} >
z^* \}$.  It follows that $\Pr \{ p \leq \udl{\bs{p}} \} \leq \Pr \{
\wh{\bs{p}}
> z^* \} \leq \f{\de}{2}$.

\epf

\sect{Proof of Theorem 7}

Theorem 7 can be shown by using the following result (a slight
modification of Hoeffding's inequality \cite{Hoeffding}) and a
similar argument as that of Theorem 3.

\beL \la{lem1Fi} Let $X_1, \cd, X_n$ be random variables with joint
distribution given by (\ref{dep}). Then, {\small $\Pr \li \{
\f{\sum_{i =1}^n X_i}{n} \geq z \ri \} \leq \exp \li ( n
\mscr{M}_{\mrm{B}} (z, p) \ri )$ } for $1 \geq z \geq p = \f{M}{N}$.
Similarly, {\small $\Pr \li \{ \f{\sum_{i =1}^n X_i}{n} \leq z \ri
\} \leq \exp \li ( n \mscr{M}_{\mrm{B}} (z, p) \ri )$ } for $0 \leq
z \leq p$. \eeL

\bpf

For $z = p$, we have {\small $\Pr \li \{ \f{\sum_{i =1}^n X_i}{n}
\geq z \ri \} \leq \exp \li ( n \mscr{M}_{\mrm{B}} (z, p) \ri ) =
1$}.  For $p < z < 1$, it was shown by Hoeffding in \cite{Hoeffding}
that {\small $\Pr \li \{ \f{\sum_{i =1}^n X_i}{n} \geq z \ri \} \leq
\exp \li ( n \mscr{M}_{\mrm{B}} (z, p) \ri )$}.  For $z = 1$,
{\small $\Pr \li \{ \f{\sum_{i =1}^n X_i}{n} \geq z \ri \} =  \Pr \{
X_i = 1, \; i = 1, \cd, n \} = \bi{M}{n} \sh \bi{N}{n} \leq p^n  =
\exp \li ( n \mscr{M}_{\mrm{B}} (1, p) \ri )$}.

For $z = 0$, {\small $\Pr \li \{ \f{\sum_{i =1}^n X_i}{n} \leq z \ri
\} = \Pr \{ X_i = 0, \; i = 1, \cd, n \} = \bi{N - M}{n} \sh
\bi{N}{n} \leq (1 - p)^n = \exp \li ( n \mscr{M}_{\mrm{B}} (0, p)
\ri )$}.   For $0 < z < p$, it was shown by Hoeffding in
\cite{Hoeffding} that {\small $\Pr \li \{ \f{\sum_{i =1}^n X_i}{n}
\leq z \ri \} \leq \exp \li ( n \mscr{M}_{\mrm{B}} (z, p) \ri )$}.
For $z = p$, we have {\small $\Pr \li \{ \f{\sum_{i =1}^n X_i}{n}
\leq z \ri \} \leq \exp \li ( n \mscr{M}_{\mrm{B}} (z, p) \ri ) =
1$}.

\epf

\sect{Proof of Theorem 8}

By the same argument as that of Theorem 1, we can show the following
lemma.

\beL \la{lemdeF}
For any $z > 0$,  {\small
\[ \Pr \{ \wh{\bs{p}} \leq z \} = \bec \Pr \{ \sum_{i = 1}^{\lc
\ga \sh z \rc - 1} X_i <
\ga \} & \tx{for} \; \ga \leq n z,\\
\Pr \li \{ \sum_{i = 1}^n X_i \leq n z \ri \} & \tx{for} \; \ga
> n z \eec \qqu  \Pr \{  \wh{\bs{p}} \geq z \} = \bec \Pr  \{ \sum_{i = 1}^{
\lf \ga \sh z \rf} X_i
\geq \ga  \} & \tx{for} \; \ga \leq n z,\\
\Pr \{ \sum_{i = 1}^{n} X_i \geq n z  \} & \tx{for} \; \ga
> n z. \eec
\]}
\eeL

By applying Lemma \ref{lemdeF}, we can show the following lemma.
\beL
 \la{Fincom383} $\Pr \{ \wh{\bs{p}} \leq  z \mid M \}$ is monotonically
decreasing with respect to $M$. Similarly, $\Pr \{ \wh{\bs{p}} \geq
z \mid M \}$ is monotonically increasing with respect to $M$.

\eeL

Now we shall introduce some new functions.  Let $p_0 < p_1 < \cd <
p_j$ be all possible values of $\wh{\bs{p}}$. Define random variable
$R$ such that $\Pr \{ R = r \} = \Pr \{ \wh{\bs{p}} = p_r \}$ for $r
= 0, 1, \cd, j$. Then, $\mscr{U}( \wh{\bs{p}} ) = \mscr{U} (p_R)$.
We denote $\mscr{U} (p_R)$ as $\mcal{U} (R)$. Clearly, $\mcal{U}
(.)$ is a non-decreasing function defined on domain $\{0, 1, \cd,
j\}$. By a linear interpolation, we can extend $\mcal{U} (.)$ as a
continuous and non-decreasing function on $[0, j]$.  Accordingly, we
can define inverse function $\mcal{U}^{-1} (.)$ such that
$\mcal{U}^{-1} (\se) = \max \{ x \in [0, j]: \mcal{U} (x ) = \se \}$
for $\mscr{U}(0) \leq \se \leq \mscr{U}(j)$. Then, $\se \geq
\mcal{U} (R) \LRA R \leq \mcal{U}^{-1} (\se) \LRA R \leq g(\se)$
where $g(\se) = \lf \mcal{U}^{-1} (\se) \rf$.

Similarly, $\mscr{L}( \wh{\bs{p}} ) = \mscr{L} (p_R)$. We denote
$\mscr{L} (p_R)$ as $\mcal{L} (R)$. Clearly, $\mcal{L} (.)$ is a
non-decreasing function defined on domain $\{0, 1, \cd, j\}$. By a
linear interpolation, we can extend $\mcal{L} (.)$ as a continuous
and non-decreasing function on $[0, j]$.  Accordingly, we can define
inverse function $\mcal{L}^{-1} (.)$ such that $\mcal{L}^{-1} (\se)
= \min \{ x \in [0, j]: \mcal{L} (x ) = \se \}$ for $\mscr{L}(0)
\leq \se \leq \mscr{L}(j)$.  Then, $\se \leq \mcal{L} (R) \LRA R
\geq \mcal{L}^{-1} (\se) \LRA R \geq h(\se)$ where $h(\se) = \lc
\mcal{L}^{-1} (\se) \rc$.

\beL \la{Flemh} Let $0 \leq r < j$.   Then, $h(m) = r + 1$ for
$\mcal{L}(r) < m \leq \mcal{L}(r + 1)$. \eeL

\bpf  Clearly, $h(m) = r + 1$ for $m = \mcal{L}(r + 1)$.  It remains
to evaluate $h(m)$ for $m$ satisfying $\mcal{L}(r) < m < \mcal{L}(r
+ 1)$.

For $m > \mcal{L}(r)$, we have $r < \mcal{L}^{-1} (m)$, otherwise $r
\geq \mcal{L}^{-1}(m)$, implying $\mcal{L}(r) \geq m$, since
$\mcal{L}(.)$ is non-decreasing and $m \notin \{\mcal{L}(r): 0 \leq
r \leq j \}$.  For $m < \mcal{L}(r + 1)$, we have $r + 1 >
\mcal{L}^{-1} (m)$, otherwise $r + 1 \leq \mcal{L}^{-1}(m)$,
implying $\mcal{L}(r + 1) \leq m$, since $\mcal{L}(.)$ is
non-decreasing and $m \notin \{\mcal{L}(r): 0 \leq r \leq j \}$.
Therefore, we have $r < \mcal{L}^{-1} (m) < r + 1$ for $\mcal{L}(r)
< m < \mcal{L}(r + 1)$. Hence, $r < \lc \mcal{L}^{-1} (m) \rc  \leq
r + 1$,  i.e., $r < h(m) \leq r + 1$. Since $h(m)$ is an integer, we
have $h(m) = r + 1$ for $\mcal{L}(r) < m < \mcal{L}(r + 1)$.

\epf

\beL \la{Flemg} Let $0 \leq r < j$.   Then, $g(m) = r$ for
$\mcal{U}(r) \leq m < \mcal{U}(r + 1)$.
 \eeL

\bpf  Clearly, $g(m) = r$ for $m = \mcal{U}(r)$.  It remains to
evaluate $g(m)$ for $m$ satisfying $\mcal{U}(r) < m < \mcal{U}(r +
1)$.

For $m > \mcal{U}(r)$, we have $r < \mcal{U}^{-1} (m)$, otherwise $r
\geq \mcal{U}^{-1}(m)$, implying $\mcal{U}(r) \geq m$, since
$\mcal{U}(.)$ is non-decreasing and $m \notin \{\mcal{U}(r): 0 \leq
r \leq j \}$.  For $m < \mcal{U}(r + 1)$, we have $r + 1 >
\mcal{U}^{-1} (m)$, otherwise $r + 1 \leq \mcal{U}^{-1}(m)$,
implying $\mcal{U}(r + 1) \leq m$, since $\mcal{U}(.)$ is
non-decreasing and $m \notin \{\mcal{U}(r): 0 \leq r \leq j \}$.
Therefore, for $\mcal{U}(r) < m < \mcal{U}(r + 1)$, we have $r <
\mcal{U}^{-1} (m)  < r + 1$. Hence, $r \leq \lf \mcal{U}^{-1} (m)
\rf  < r + 1$, i.e., $r \leq g(m) < r + 1$. Since $g(m)$ is an
integer, we have $g(m) = r$ for $\mcal{U}(r) <  m < \mcal{U}(r +
1)$.

\epf

Noting that $\Pr \{ M \geq \mscr{U}( \wh{\bs{p}} ) \mid M \} = \Pr
\{ M \geq \mcal{U}( R ) \mid M \}$, we have $\Pr \{ M \geq \mscr{U}(
\wh{\bs{p}} ) \mid M \}  = \Pr \{ R \leq g(M) \mid M \}$.   Let $0
\leq r < j$. By Lemma \ref{Flemg}, we have that $g(m) = r$ for
$\mcal{U} (r) \leq m < \mcal{U} (r + 1)$.  Observing that $\Pr \{ M
\geq \mscr{U}( \wh{\bs{p}} ) \mid M \}  = 0$ for $0 \leq M <
\mscr{U} (0)$ and that $\Pr \{ M \geq \mscr{U}( \wh{\bs{p}} ) \mid M
\}  = 1$ for $\mscr{U} (j) \leq M \leq N$, we have that the maximum
of $\Pr \{ M \geq \mscr{U}( \wh{\bs{p}} ) \mid M \}$ with respect to
$M \in [a, b]$ is achieved on $\bigcup_{r = 0}^{j - 1} \{ m \in [a,
b] : \mcal{U} (r) \leq m \leq \mcal{U} (r + 1) \} \cup \{ a, b \}$.
Now consider the range $\{ m \in [a, b] : \mcal{U} (r) \leq m \leq
\mcal{U} (r + 1) \}$ of $M$.  We only consider the non-trivial
situation that $\mcal{U} (r) < \mcal{U} (r + 1)$.  For $\mcal{U} (r)
\leq M < \mcal{U} (r + 1)$, we have
\[
\Pr \{ M \geq \mscr{U}( \wh{\bs{p}} ) \mid M \} = \Pr \{ R \leq g(M)
\mid M \} = \Pr \{ R \leq r \mid M \} = \Pr \{ \wh{\bs{p}} \leq p_r
\mid M \},
\]
which is non-increasing for this range of $M$ as can be seen from
Lemma \ref{Fincom383}. By virtue of such monotonicity, we can
characterize the maximizer of $\Pr \{ M \geq \mscr{U}( \wh{\bs{p}} )
\mid M \}$ with respect to $M$ on the set $\{ m \in [a, b] :
\mcal{U} (r) \leq m \leq \mcal{U} (r + 1) \}$ as follows.

Case (i): $b < \mcal{U} (r)$ or $a > \mcal{U} (r + 1)$.  This is
trivial.

Case (ii): $a < \mcal{U} (r) \leq b \leq \mcal{U} (r + 1)$. The
maximizer must be among $\{\mcal{U} (r), \; b \}$.

Case (iii): $\mcal{U} (r) \leq a \leq  b \leq \mcal{U} (r + 1)$. The
maximizer must be among $\{a, \; b \}$.

Case (iv): $\mcal{U} (r)  \leq  a \leq \mcal{U} (r + 1) < b$. The
maximizer must be among $\{a, \; \mcal{U} (r + 1) \}$.

Case (v): $a < \mcal{U} (r)  \leq \mcal{U} (r + 1) < b$. The
maximizer must be among $\{\mcal{U} (r),  \; \mcal{U} (r + 1) \}$.

In summary, the maximizer must be among $\{\mcal{U} (r),  \;
\mcal{U} (r + 1), a, b \} \cap [a, b]$. It follows that the
statement on $\Pr \{ M \geq \mscr{U}( \wh{\bs{p}} ) \mid M \}$ is
established.

\bsk

Next, we consider $\Pr \{ M \leq \mscr{L}( \wh{\bs{p}} ) \mid M \}$.
Noting that $\Pr \{ M \leq \mscr{L}( \wh{\bs{p}} ) \mid M \} = \Pr
\{ M \leq \mcal{L}( R ) \mid M \}$, we have $\Pr \{ M \leq \mscr{L}(
\wh{\bs{p}} ) \mid M \}  = \Pr \{ R \geq h(M) \mid M \}$. Let $0
\leq r < j$. By Lemma \ref{Flemh}, we have that $h(m) = r + 1$ for
$\mcal{L} (r) < m \leq \mcal{L} (r + 1)$.  Observing that $\Pr \{ M
\leq \mscr{L}( \wh{\bs{p}} ) \mid M \} = 1$ for $0 \leq M \leq
\mscr{L} (0)$ and that $\Pr \{ M \leq \mscr{L}( \wh{\bs{p}} ) \mid M
\} = 0$ for $\mscr{L} (j) < M \leq N$, we have that the maximum of
$\Pr \{ M \leq \mscr{L}( \wh{\bs{p}} ) \mid M \}$ with respect to $M
\in [a, b]$ is achieved on $\bigcup_{r = 0}^{j - 1} \{ m \in [a, b]
: \mcal{L} (r) \leq m \leq \mcal{L} (r + 1) \} \cup \{ a, b \}$. Now
consider the range $\{ m \in [a, b] : \mcal{L} (r) \leq m \leq
\mcal{L} (r + 1) \}$ of $M$.  We only consider the non-trivial
situation that $\mcal{L} (r) < \mcal{L} (r + 1)$.  For $\mcal{L} (r)
< M \leq \mcal{L} (r + 1)$, we have
\[
\Pr \{ M \leq \mscr{L}( \wh{\bs{p}} ) \mid M \} = \Pr \{ R \geq h(M)
\mid M \} = \Pr \{ R \geq r + 1 \mid M \} = \Pr \{ \wh{\bs{p}} \geq
p_{r+ 1} \mid M \},
\]
which is non-decreasing for this range of $M$ as can be seen from
Lemma \ref{Fincom383}. By virtue of such monotonicity, we can
characterize the maximizer of $\Pr \{ M \leq \mscr{L}( \wh{\bs{p}} )
\mid M \}$ with respect to $M$ on the set $\{ m \in [a, b] :
\mcal{L} (r) \leq m \leq \mcal{L} (r + 1) \}$ as follows.

Case (i): $b < \mcal{L} (r)$ or $a > \mcal{L} (r + 1)$.  This is
trivial.

Case (ii): $a < \mcal{L} (r) \leq b \leq \mcal{L} (r + 1)$. The
maximizer must be among $\{\mcal{L} (r), \; b \}$.

Case (iii): $\mcal{L} (r) \leq a \leq  b \leq \mcal{L} (r + 1)$. The
maximizer must be among $\{a, \; b \}$.

Case (iv): $\mcal{L} (r)  \leq  a \leq \mcal{L} (r + 1) < b$. The
maximizer must be among $\{a, \; \mcal{L} (r + 1) \}$.

Case (v): $a < \mcal{L} (r)  \leq \mcal{L} (r + 1) < b$. The
maximizer must be among $\{\mcal{L} (r),  \; \mcal{L} (r + 1) \}$.

In summary, the maximizer must be among $\{\mcal{L} (r),  \;
\mcal{L} (r + 1), a, b \} \cap [a, b]$. It follows that the
statement on $\Pr \{ M \leq \mscr{L}( \wh{\bs{p}} ) \mid M \}$ is
established.

This concludes the proof of Theorem \ref{THM3}.

\sect{Proof of Theorem 10}

The theorem can be established by showing the following lemmas.

\beL $\Pr \{ M > \bs{M}_u \} \leq \f{\de}{2}$. \eeL

\bpf

Since $X_i$ must be either $0$ or $1$ and $\ga$ is an integer, we
have {\small $\Pr \{ \sum_{i = 1}^{\lc \ga \sh z \rc - 1} X_i < \ga
\} \leq \Pr \{ \sum_{i = 1}^{\lc \ga \sh z \rc} X_i \leq \ga \} $}.
Hence, by Lemma \ref{lemdeF}, {\small
\[ \Pr \{ \wh{\bs{p}} \leq z \} \leq \bec \Pr \{ \sum_{i = 1}^{\lc
\ga \sh z \rc} X_i \leq
\ga \} & \tx{for} \; \ga \leq n z,\\
\Pr \li \{ \sum_{i = 1}^n X_i \leq n z \ri \} & \tx{for} \; \ga
> n z \eec
\]}
and $\Pr \{  \wh{\bs{p}} \leq z \} \leq \mscr{G}( z, p)$, where
\[
\mscr{G}( z, p) = \bec \sum_{i = 0}^{\ga} \bi{M}{i} \bi{N - M}{\lc
\ga \sh z \rc - i}
\sh \bi{N}{\lc \ga \sh z \rc}  &  \tx{for} \; \f{\ga}{n} \leq z \leq 1,\\
\sum_{i = 0}^{\lf n z \rf} \bi{M}{i} \bi{N - M}{n - i} \sh \bi{N}{n}
&  \tx{for} \; 0 \leq z < \f{\ga}{n} \eec
\]
with $p = \f{M}{N}$.  Let $z^* \in [0, 1]$ be the largest number
such that $\Pr \{ \wh{\bs{p}} < z^* \} \leq \f{\de}{2}$.  Since
$\wh{\bs{p}}$ is a discrete random variable bounded in $[0, 1]$, it
must be true that $\Pr \{ \wh{\bs{p}} \leq z^* \} > \f{\de}{2}$.
Observing that $\mscr{G} (z, p)$ is monotonically decreasing with
respect to $p \in \{\f{i}{N} : i = 0, 1, \cd, N \}$, we have {\small
\[ \{ p \geq \ovl{\bs{p}} \} \subseteq  \li \{ \mscr{G} (\wh{\bs{p}}, p) \leq \mscr{G}
(\wh{\bs{p}}, \ovl{\bs{p}}) \leq \f{\de}{2} \ri \} \subseteq \li \{
\mscr{G} (\wh{\bs{p}}, p) \leq \f{\de}{2} \ri \}
\]}
where $\ovl{\bs{p}} = \f{\bs{M}_u + 1 } {N}$. Noting that
$\f{\de}{2} < \Pr \{ \wh{\bs{p}} \leq z^* \} \leq \mscr{G} (z^*, p)$
and that $\mscr{G} (z, p)$ is non-decreasing with respect to $z \in
(0, 1)$, we have $ \{ p \geq \ovl{\bs{p}} \} \subseteq \{ \mscr{G}
(\wh{\bs{p}}, p) \leq \f{\de}{2} \} \subseteq \{ \mscr{G}
(\wh{\bs{p}}, p)  < \mscr{G} (z^*, p) \} \subseteq \{ \wh{\bs{p}} <
z^* \}$.  It follows that $\Pr \{ p \geq \ovl{\bs{p}} \} \leq \Pr \{
\wh{\bs{p}} < z^* \} \leq \f{\de}{2}$, which implies that $\Pr \{ M
> \bs{M}_u \} \leq \f{\de}{2}$.

\epf

\beL $\Pr \{ M < \bs{M}_l \} \leq \f{\de}{2}$. \eeL

\bpf  By Lemma \ref{lemdeF}, we have $\Pr \{  \wh{\bs{p}} \geq z \}
= \mscr{H}( z, p)$, where
\[
\mscr{H}( z, p) = \bec \sum_{i = \ga}^{\lf \ga \sh z \rf}
\bi{M}{i} \bi{N - M}{\lf \ga \sh z \rf - i} \sh \bi{N}{\lf \ga \sh z \rf}
&  \tx{for} \; \f{\ga}{n} \leq z \leq 1,\\
\sum_{i = \lc n z \rc}^n  \bi{M}{i} \bi{N - M}{n - i} \sh \bi{N}{n}
&  \tx{for} \; 0 \leq z < \f{\ga}{n} \eec
\]
with $p = \f{M}{N}$.  Let $z^* \in [0, 1]$ be the smallest number
such that $\Pr \{ \wh{\bs{p}} > z^* \} \leq \f{\de}{2}$.  Since
$\wh{\bs{p}}$ is a discrete random variable bounded in $[0, 1]$, it
must be true that $\Pr \{ \wh{\bs{p}} \geq z^* \} > \f{\de}{2}$.
Observing that $\mscr{H} (z, p)$ is monotonically increasing with
respect to $p \in \{\f{i}{N} : i = 0, 1, \cd, N \}$, we have {\small
\[ \{ p \leq \udl{\bs{p}} \} \subseteq  \li \{ \mscr{H} (\wh{\bs{p}}, p) \leq \mscr{H}
(\wh{\bs{p}}, \udl{\bs{p}}) \leq \f{\de}{2} \ri \} \subseteq \li \{
\mscr{H} (\wh{\bs{p}}, p) \leq \f{\de}{2} \ri \}
\]}
where $\udl{\bs{p}} = \f{\bs{M}_l - 1 } {N}$.  Noting that
$\f{\de}{2} < \Pr \{ \wh{\bs{p}} \geq z^* \} = \mscr{H} (z^*, p)$
and that $\mscr{H} (z, p)$ is non-increasing with respect to $z \in
(0, 1)$, we have $ \{ p \leq \udl{\bs{p}} \} \subseteq \{ \mscr{H}
(\wh{\bs{p}}, p) \leq \f{\de}{2} \} \subseteq \{ \mscr{H}
(\wh{\bs{p}}, p)  < \mscr{H} (z^*, p) \} \subseteq \{ \wh{\bs{p}} >
z^* \}$.  It follows that $\Pr \{ p \leq \udl{\bs{p}} \} \leq \Pr \{
\wh{\bs{p}}
> z^* \} \leq \f{\de}{2}$, which implies $\Pr \{ M < \bs{M}_l \} \leq \f{\de}{2}$.

\epf

\sect{Proof of Theorem 11}

In the case of $M < \ga$, we have $\mbf{n} = n$ and $\f{\ga}{p} >
\f{\ga N}{M} > N$, from which the theorem immediately follows.  It
remains to show the theorem for the case of $M \geq \ga$.   Notice
that \bee \bb{E} [ \mbf{n} ] & = & n \Pr \li \{ \sum_{i = 1}^n X_i <
\ga \ri \}  + \sum_{m = 1}^n m \Pr \li \{ \sum_{i = 1}^{m - 1} X_i <
\ga, \; \sum_{i = 1}^m X_i \geq \ga \ri \}\\
& < & n \Pr \li \{ \sum_{i = 1}^n X_i < \ga \ri \}  + \sum_{m = 1}^n
n \Pr \li \{ \sum_{i = 1}^{m - 1} X_i < \ga, \; \sum_{i = 1}^m X_i
\geq \ga \ri \} = n.  \eee  Since  $M = \sum_{i = 1}^N X_i \geq \ga$
and $X_i$ is non-negative, we have {\small
\[ \bigcup_{m = n + 1}^N \li \{ \sum_{i = 1}^{m - 1} X_i < \ga, \;
\sum_{i = 1}^m X_i \geq \ga \ri \}   = \li \{ \sum_{i = 1}^n X_i <
\ga \ri \}.
\]}
Hence, {\small \bee \bb{E} [ \mbf{n} ] & = &  \sum_{m = n + 1}^N n
\Pr \li \{ \sum_{i = 1}^{m - 1} X_i < \ga, \; \sum_{i = 1}^m X_i
\geq \ga \ri \} + \sum_{m = 1}^n m \Pr \li \{ \sum_{i = 1}^{m - 1}
X_i <
\ga, \; \sum_{i = 1}^m X_i \geq \ga \ri \}\\
& < & \sum_{m = 1}^N m \Pr \li \{ \sum_{i = 1}^{m - 1} X_i < \ga, \;
\sum_{i = 1}^m X_i \geq \ga \ri \} = \bb{E} [ \mbf{m} ], \eee} where
$\mbf{m}$ is the sample number of the classical {\it inverse
sampling} scheme with the following stopping rule: Sampling without
replacement is continued until $\ga$ units possessing the attribute
have been observed. By the definition of the classical inverse
sampling, we have $\sum_{i = 1}^{\mbf{m}} X_i = \ga$. Noting that
$X_1, X_2, \cd, X_N$ are identical but dependent Bernoulli random
variables with common mean $p$ and that $\{ \mbf{n} \geq k \}$
depends only on $X_1, \cd, X_{k - 1}$ for $1 \leq k \leq N$, we can
conclude that Wald's equation still applies. Hence, $\bb{E} [
\sum_{i = 1}^{\mbf{m} } X_i ] = \bb{E} [ \mbf{m} ] \; \bb{E} [ X_i ]
= \ga$, which implies that $\bb{E} [ \mbf{m} ] = \f{\ga}{\bb{E} [
X_i ]} = \f{\ga}{p}$. Since $\bb{E} [ \mbf{n} ]$ is less than both
$n$ and $\bb{E} [ \mbf{m} ]$ as shown above, we have {\small $\bb{E}
[ \mbf{n} ] < \min \{ n, \f{\ga}{p} \}$}.  This completes the proof
of the theorem.

\sect{Proof of Theorem 12}

The theorem can be established by showing the following lemmas.

\beL $\Pr \{ \lm \geq \ovl{\bs{\lm}} \} \leq \f{\de}{2}$. \eeL

\bpf

By Theorem 1, we have {\small \[ \Pr \{ \wh{\bs{\lm}} \leq z \} =
\bec \Pr \{ \sum_{i = 1}^{\lc \ga \sh z \rc - 1} X_i <
\ga \} & \tx{for} \; \ga \leq n z,\\
\Pr \li \{ \sum_{i = 1}^n X_i \leq n z \ri \} & \tx{for} \; \ga
> n z. \eec
\]}
and thus $\Pr \{  \wh{\bs{\lm}} \leq z \} = \mscr{G}( z, \lm)$,
where
\[
\mscr{G}( z, \lm) = \bec \sum_{i = 0}^{\ga - 1} \f{1}{ i! } [( \lc
\f{\ga}{z} \rc - 1) \lm ]^i \exp ( - ( \lc \f{\ga}{z} \rc
- 1) \lm )  &  \tx{for} \; z \geq \f{\ga}{n},\\
\sum_{i = 0}^{\lf n z \rf} \f{1}{ i! } ( n \lm )^i \exp ( - n \lm )
&  \tx{for} \; 0 \leq z < \f{\ga}{n}. \eec
\]
Let $z^* \geq 0$ be the largest number such that $\Pr \{
\wh{\bs{\lm}} < z^* \} \leq \f{\de}{2}$.  Since $\wh{\bs{\lm}}$ is a
non-negative discrete random variable, it must be true that $\Pr \{
\wh{\bs{\lm}} \leq z^* \} > \f{\de}{2}$. Observing that $\mscr{G}
(z, \lm)$ is monotonically decreasing with respect to $\lm \in (0,
\iy)$, we have {\small \[ \{ \lm \geq \ovl{\bs{\lm}} \} \subseteq
\li \{ \mscr{G} (\wh{\bs{\lm}}, \lm) \leq \mscr{G} (\wh{\bs{\lm}},
\ovl{\bs{\lm}}) = \f{\de}{2} \ri \} \subseteq \li \{ \mscr{G}
(\wh{\bs{\lm}}, \lm) \leq \f{\de}{2} \ri \}.
\]}
Noting that $\f{\de}{2} < \Pr \{ \wh{\bs{\lm}} \leq z^* \} =
\mscr{G} (z^*, \lm)$ and that $\mscr{G} (z, \lm)$ is non-decreasing
with respect to $z \in (0, \iy)$, we have $ \{ \lm \geq
\ovl{\bs{\lm}} \} \subseteq \{ \mscr{G} (\wh{\bs{\lm}}, \lm) \leq
\f{\de}{2} \} \subseteq \{ \mscr{G} (\wh{\bs{\lm}}, \lm) < \mscr{G}
(z^*, \lm) \} \subseteq \{ \wh{\bs{\lm}} < z^* \}$. It follows that
$\Pr \{ \lm \geq \ovl{\bs{\lm}} \} \leq \Pr \{ \wh{\bs{\lm}} < z^*
\} \leq \f{\de}{2}$.

\epf

\beL $\Pr \{ \lm \leq \udl{\bs{\lm}} \} \leq \f{\de}{2}$. \eeL

\bpf

 By Theorem 1, we have {\small
\[
\Pr \{  \wh{\bs{\lm}} \geq z \} = \bec \Pr  \{ \sum_{i = 1}^{ \lf
\ga \sh z \rf} X_i
\geq \ga  \} & \tx{for} \; \ga \leq n z,\\
\Pr \{ \sum_{i = 1}^{n} X_i \geq n z  \} & \tx{for} \; \ga
> n z \eec
\]}
and thus $\Pr \{  \wh{\bs{\lm}} \geq z \} = \mscr{H}( z, \lm)$ where
\[
\mscr{H}( z, \lm) = \bec \sum_{i = \ga}^{\iy} \f{1}{ i! } ( \lf
\f{\ga}{z} \rf  \lm )^i \exp ( - \lf \f{\ga}{z} \rf \lm
)  &  \tx{for} \; z \geq \f{\ga}{n},\\
\sum_{i = \lc n z \rc}^\iy \f{1}{ i! } ( n \lm )^i \exp ( - n \lm )
&  \tx{for} \; 0 \leq z < \f{\ga}{n}. \eec
\]
Let $z^* \geq 0$ be the smallest number such that $\Pr \{
\wh{\bs{\lm}} > z^* \} \leq \f{\de}{2}$.  Since $\wh{\bs{\lm}}$ is a
non-negative discrete random variable, it must be true that $\Pr \{
\wh{\bs{\lm}} \geq z^* \} > \f{\de}{2}$.  Observing that $\mscr{H}
(z, \lm)$ is monotonically increasing with respect to $\lm \in (0,
\iy)$, we have {\small \[ \{ \lm \leq \udl{\bs{\lm}} \} = \{ \lm
\leq \udl{\bs{\lm}}, \; \mbf{k} > 0 \} \subseteq  \li \{ \mscr{H}
(\wh{\bs{\lm}}, \lm) \leq \mscr{H} (\wh{\bs{\lm}}, \udl{\bs{\lm}}) =
\f{\de}{2} \ri \} \subseteq \li \{ \mscr{H} (\wh{\bs{\lm}}, \lm)
\leq \f{\de}{2} \ri \}.
\]}
Noting that $\f{\de}{2} < \Pr \{ \wh{\bs{\lm}} \geq z^* \} =
\mscr{H} (z^*, \lm)$ and that $\mscr{H} (z, \lm)$ is non-increasing
with respect to $z \in (0, \iy)$, we have $\{ \lm \leq
\udl{\bs{\lm}} \} \subseteq \{ \mscr{H} (\wh{\bs{\lm}}, \lm) \leq
\f{\de}{2} \} \subseteq \{ \mscr{H} (\wh{\bs{\lm}}, \lm) < \mscr{H}
(z^*, \lm) \} \subseteq \{ \wh{\bs{\lm}} > z^* \}$.  It follows that
$\Pr \{ \lm \leq \udl{\bs{\lm}} \} \leq \Pr \{ \wh{\bs{\lm}}
> z^* \} \leq \f{\de}{2}$.

\epf

\sect{Proof of Theorem 13}

By the same method as that of Lemma \ref{lem11}, we have \beL
\la{lem11Fi} $\Pr \{ \wh{\bs{\mu}} \geq (1 + \vep_r) \mu \} <
\f{\de}{2}$ for any $\mu \in \li ( p^\star, 1 \ri )$. \eeL

\beL \la{lem24} Let $0 < \vep_r < 1, \; z = (1 - \vep_r) \mu$ and
$\vep^\prime = 1 - \f{\ga (1 - \vep_r)}{\ga + \vep_r - 1}$. Suppose
$\ga
> \f{ 1 - \vep_r } {\vep_r }$. Then, {\small $\mscr{M}_{\mrm{I}} \li (
\f{z \ga}{ \ga - z }, \mu \ri ) < \mscr{M}_{\mrm{I}} \li ( (1 -
\vep^\prime) \mu, \mu \ri )$} for any $\mu \in (0, 1)$. \eeL

\bpf

As a consequence of $\ga > \f{ 1 - \vep_r } {\vep_r }$, we have $0 <
\f{z \ga}{ \ga - z } < \mu$ for any $\mu \in (0, 1)$.
 Since $ \mscr{M}_{\mrm{I}} (w, \mu)$ is monotonically increasing with respect to
$w \in (0, \mu)$, it suffices to show that $\f{z \ga}{ \ga  - z } <
(1 - \vep^\prime) \mu$ for any $\mu \in (0, 1)$. That is, to show
$\f{(1 - \vep_r) \mu \ga}{ \ga - (1 - \vep_r) \mu } < (1 -
\vep^\prime) \mu, \; \fa \mu \in (0, 1)$, i.e., $\f{(1 - \vep_r)
\ga}{ \ga  - (1 - \vep_r) \mu } < 1 - \vep^\prime, \; \fa \mu \in
(0, 1)$. This follows from the definition of $\vep^\prime$.  \epf

\beL  \la{lem12Fi}  $\Pr \{ \wh{\bs{\mu}} \leq (1 - \vep_r) \mu \} <
\f{\de}{2}$ for any $\mu \in \li ( p^\star, 1 \ri )$. \eeL

\bpf To prove the lemma, we shall consider the following two cases:

Case (i): $(1 - \vep_r) \mu \geq \f{\ga}{n}$;

Case (ii): $(1 - \vep_r) \mu < \f{\ga}{n}$.

\bsk

For Case (i),  applying Theorem 1 with $z = (1 - \vep_r) \mu \geq
\f{\ga}{n}$, we have {\small $\Pr \{ \wh{\bs{\mu}} \leq (1 - \vep_r)
\mu \}  =  \Pr \li \{ \sum_{i = 1}^{\lc \ga \sh z \rc - 1} X_i < \ga
\ri \}$}. As a consequence of $\ga > \f{ 1 - \vep_r } {\vep_r }$, we
have $0 < \f{\ga}{\lc \ga \sh z \rc - 1} \leq \f{z \ga}{ \ga - z } <
\mu$ for any $\mu \in (0, 1)$.  Hence, applying Lemma \ref{lem1}, we
have {\small \bel \Pr \{ \wh{\bs{\mu}} \leq (1 - \vep_r) \mu \} &
\leq  & \exp \li (  ( \lc \ga \sh z \rc - 1) \; \mscr{M}_{\mrm{B}}
\li ( \f{\ga}{\lc \ga \sh z \rc - 1}, \mu
 \ri ) \ri ) \nonumber\\
 &  = & \exp \li (  \ga \; \mscr{M}_{\mrm{I}} \li (
\f{\ga}{\lc \ga \sh z \rc - 1}, \mu
 \ri ) \ri ) \leq \exp \li (  \ga \; \mscr{M}_{\mrm{I}} \li ( \f{\ga z}{\ga - z}, \mu
 \ri ) \ri ). \nonumber \eel}
By Lemmas \ref{lem24} and \ref{lem3}, \bee \Pr \{ \wh{\bs{\mu}} \leq
(1 - \vep_r) \mu \} & \leq & \exp \li ( \ga \mscr{M}_{\mrm{I}} ( (1
- \vep^\prime) \mu, \mu ) \ri ) \leq  \exp \li ( \ga
\mscr{M}_{\mrm{I}} ( (1 - \vep^\prime) p^\star, p^\star ) \ri )\\
&  = & \exp \li ( \ga \mscr{M}_{\mrm{I}} \li ( \f{\ga (1 - \vep_r)
p^\star}{\ga + \vep_r - 1}, p^\star \ri ) \ri ) < \f{\de}{2}.  \eee

For Case (ii), applying Theorem 1 with $z = (1 - \vep_r) \mu <
\f{\ga}{n}$ and by a similar argument as that of Lemma \ref{lem12},
we have $\Pr \{ \wh{\bs{\mu}} \leq (1 - \vep_r) \mu \} <
\f{\de}{2}$.

In summary, we have shown $\Pr \{ \wh{\bs{\mu}} \leq (1 - \vep_r)
\mu \} < \f{\de}{2}$ for all cases.  The lemma is thus proved.

\epf

By the same method as that of Lemma \ref{lem13}, we have

\beL  \la{lem13Fi}  $\Pr \{ \wh{\bs{\mu}} \geq \mu  + \vep_a \} <
\f{\de}{2}$ for any $\mu \in \li (0, p^\star \ri ]$. \eeL

\beL \la{lem26} Let $z = \mu - \vep$.  Suppose $0 < \f{z \ga}{ \ga
- z } < \mu$.  Then,  $\mscr{M}_{\mrm{I}} \li ( \f{z \ga}{ \ga  - z
}, \mu \ri )$ is monotonically increasing with respect to $\mu \in
\li ( \vep, \f{1}{2} \ri )$.   \eeL

\bpf

Note that {\small $\f{ \pa \mscr{M}_{\mrm{I}} (w, \mu) }{ \pa \mu} =
\f{1}{\mu}  - \f{1}{w^2} \ln \li ( \f{ 1 - \mu } { 1 - w } \ri )
\f{\pa w}{\pa \mu} - \li ( \f{1}{w} - 1 \ri ) \f{1}{1 - \mu}$},
where $w = \f{z \ga}{ \ga - z } = - \ga + \f{\ga^2}{ \ga  - z }$ and
$\f{\pa w}{\pa \mu} = \f{\ga^2}{ (\ga - z )^2 }  = \f{w^2}{ z^2}$.
Hence, {\small \bee \f{ \pa \mscr{M}_{\mrm{I}} (w, \mu) }{ \pa \mu}
& = & \f{1}{\mu} + \f{1}{ z^2} \ln \li ( \f{ 1 - w } { 1 - \mu } \ri
) - \li ( \f{1}{w} - 1 \ri ) \f{1}{1 - \mu} >  \f{1}{\mu}  + \f{1}{
z^2}
\li ( \f{ \mu - w } { 1 - w } \ri )  - \li ( \f{1}{w} - 1 \ri ) \f{1}{1 - \mu}\\
& = &  \f{1}{ z^2}  \li ( \f{ \mu - w } { 1 - w } \ri ) - \f{\mu -
w}{ \mu (1 - \mu) w} > 0 \eee} if $z^2 ( 1 - w ) < \mu (1 - \mu) w$,
i.e., {\small $z^2 \li ( \f{1}{w} - 1 \ri ) < \mu (1 - \mu) \LRA z^2
\li ( \f{1}{z} - \f{1}{\ga} - 1 \ri ) < \mu (1 - \mu) \LRA z (1 - z)
-  \f{z^2}{\ga} <  \mu (1 - \mu)$} since $\f{1}{w} = \f{1}{z} -
\f{1}{\ga}$.

\epf

\beL \la{lem14Fi} $\Pr \{ \wh{\bs{\mu}} \leq \mu - \vep_a \} <
\f{\de}{2}$ for any $\mu \in \li (0,  p^\star \ri ]$.

\eeL

\bpf

To prove the lemma, we shall consider the following three cases:

Case (i): $\mu < \vep_a$;

Case (ii): $\mu - \vep_a \geq \f{\ga}{n}$;

Case (iii): $0 \leq \mu - \vep_a < \f{\ga}{n}$.

\bsk

For Case (i), it is evident that $\Pr \{ \wh{\bs{\mu}} \leq \mu -
\vep_a \} = 0 < \f{\de}{2}$.

For Case (ii), applying Theorem 1 with $z = \mu - \vep_a \geq
\f{\ga}{n}$, we have {\small $\Pr \{ \wh{\bs{\mu}} \leq \mu - \vep_a
\}  = \Pr \li \{ \sum_{i = 1}^{\lc \ga \sh z \rc - 1} X_i < \ga \ri
\}$}.  From the definition of $z$ and the assumption that $\ga >
\f{1 - \vep_r}{\vep_r}$, we see that $0 < \f{\ga}{\lc \ga \sh z \rc
- 1} \leq \f{z \ga}{ \ga - z } < \mu$ for any $\mu \in (0,
p^\star]$.  Hence, it follows from  Lemma \ref{lem1} that  \bel \Pr
\{ \wh{\bs{\mu}} \leq \mu - \vep_a \} &  \leq & \exp \li (  ( \lc
\ga \sh z \rc - 1) \; \mscr{M}_{\mrm{B}} \li ( \f{\ga}{\lc \ga \sh z
\rc - 1},
\mu \ri ) \ri ) \nonumber\\
 & = & \exp \li (  \ga \; \mscr{M}_{\mrm{I}} \li (
\f{\ga}{\lc \ga \sh z \rc - 1}, \mu
 \ri ) \ri ) \leq \exp \li (  \ga \; \mscr{M}_{\mrm{I}} \li ( \f{\ga z}{\ga - z}, \mu
 \ri ) \ri ). \nonumber \eel  Invoking Lemma \ref{lem26},
 we have {\small \bee \Pr \{ \wh{\bs{\mu}} \leq
\mu - \vep_a \} & \leq & \exp \li (  \ga \; \mscr{M}_{\mrm{I}} \li (
\f{\ga (p^\star - \vep_a) }{\ga - (p^\star - \vep_a) }, p^\star \ri
) \ri ) = \exp \li (  \ga \; \mscr{M}_{\mrm{I}} \li ( \f{\ga (1 -
\vep_r) p^\star }{\ga - (1 - \vep_r) p^\star }, p^\star \ri ) \ri )\\
& < & \exp \li (  \ga \; \mscr{M}_{\mrm{I}} ( ( 1 - \vep^\prime)
p^\star, p^\star ) \ri ) =  \exp \li ( \ga \mscr{M}_{\mrm{I}}  \li (
\f{\ga (1 - \vep_r) p^\star}{\ga + \vep_r - 1}, p^\star \ri ) \ri )
< \f{\de}{2} \eee} where $\vep^\prime$ is defined in Lemma
\ref{lem24}.

For Case (iii), applying Theorem 1 with $z = \mu - \vep_a <
\f{\ga}{n}$ and by an argument as that of Lemma \ref{lem14}, we have
$\Pr \{ \wh{\bs{\mu}} \leq \mu - \vep_a \} < \f{\de}{2}$. In
summary, we have shown $\Pr \{ \wh{\bs{\mu}} \leq \mu - \vep_a \} <
\f{\de}{2}$ for all cases. The lemma is thus proved.

 \epf

\bsk

Finally, the proof of Theorem 13 can be accomplished by a similar
argument as that of Theorem 3.

\sect{Proof of Theorem 14}

The theorem can be established by showing the following lemmas.

\beL $\Pr \{ \mu \geq \ovl{\bs{\mu}} \} \leq \f{\de}{2}$. \eeL

\bpf

For $\mu > z \geq \f{\ga}{n}$, by Theorem 1 and Lemma \ref{lem1}, we
have {\small \bee \Pr \{ \wh{\bs{\mu}} \leq z \} & = & \Pr \li \{
\sum_{i = 1}^{\li \lc \ga \sh z \ri \rc - 1}  X_i < \ga \ri \} \leq
\exp \li ( (\li \lc \ga \sh z \ri \rc - 1) \mscr{M}_{\mrm{B}} \li (
\f{\ga}{\li \lc \ga \sh z \ri \rc  - 1},
\mu \ri ) \ri )\\
&  = & \exp \li (\ga \mscr{M}_{\mrm{I}} \li ( \f{\ga}{\li \lc \ga
\sh z \ri \rc  - 1}, \mu \ri ) \ri ) \leq \exp \li ( \ga
\mscr{M}_{\mrm{I}} \li ( \f{z \ga}{ \ga - z }, \mu \ri ) \ri ) \eee}
where the last inequality is due to $\f{\ga}{\li \lc \ga \sh z \ri
\rc  - 1} \leq \f{z \ga}{ \ga - z }$ and
 the fact that $\mscr{M}_{\mrm{I}} (z, \mu)$  is monotonically increasing with respect
 to $z \in (0, \mu)$.
 For $\mu > z$ and $0 \leq z <  \f{\ga}{n}$, by Theorem 1 and Lemma \ref{lem1},
  we have {\small $\Pr \li \{  \wh{\bs{\mu}} \leq z \ri \} = \Pr \li
\{  \f{ \sum_{i = 1}^{n} X_i  } { n } \leq z \ri \} \leq \exp( n
\mscr{M}_{\mrm{B}} ( z, \mu) )$}.  Therefore, $\Pr \{ \wh{\bs{\mu}}
\leq z \} \leq \mscr{G}( z, \mu)$, where
\[
\mscr{G}( z, \mu) = \bec \exp \li ( \ga \mscr{M}_{\mrm{I}} \li (
\f{z
\ga}{ \ga - z }, \mu \ri ) \ri )   &  \tx{for} \; \f{\ga}{n} \leq z < \mu,\\
\exp( n \mscr{M}_{\mrm{B}} ( z, \mu) )  & \tx{for} \; 0 \leq z <
\f{\ga}{n}, \; z < \mu. \eec
\]
Let $z^* \in [0, 1]$ be the largest number such that $\Pr \{
\wh{\bs{\mu}} < z^* \} \leq \f{\de}{2}$.  Then, it must be true that
either $\Pr \{ \wh{\bs{\mu}} \leq z^* \} > \f{\de}{2}$ or $\Pr \{
\wh{\bs{\mu}} \leq z^* \} = \f{\de}{2}$. Observing that $\mscr{G}
(z, \mu)$ is monotonically decreasing with respect to $\mu \in (z,
1)$, we have {\small \[ \{ \mu \geq \ovl{\bs{\mu}} \} = \{ \mu \geq
\ovl{\bs{\mu}} \geq \wh{\bs{\mu}}, \; \mbf{k} < \mbf{n} \} \subseteq
\li \{ \mscr{G} (\wh{\bs{\mu}}, \mu) \leq \mscr{G} (\wh{\bs{\mu}},
\ovl{\bs{\mu}}) = \f{\de}{2}, \; \mu \geq \ovl{\bs{\mu}} \geq
\wh{\bs{\mu}} \ri \} \subseteq \li \{ \mscr{G} (\wh{\bs{\mu}}, \mu)
\leq \f{\de}{2}, \; \wh{\bs{\mu}} \leq \mu \ri \}.
\]}
In the case of $\Pr \{ \wh{\bs{\mu}} \leq z^* \} > \f{\de}{2}$, we
have $\f{\de}{2} < \Pr \{ \wh{\bs{\mu}} \leq z^* \} \leq \mscr{G}
(z^*, \mu)$.  Since $\mscr{G} (z, \mu)$ is increasing with respect
to $z \in (0, \mu)$, we have $\{ \mu \geq \ovl{\bs{\mu}} \}
\subseteq \{ \mscr{G} (\wh{\bs{\mu}}, \mu) \leq \f{\de}{2}, \;
\wh{\bs{\mu}} \leq \mu  \} \subseteq \{ \mscr{G} (\wh{\bs{\mu}},
\mu) < \mscr{G} (z^*, \mu), \; \wh{\bs{\mu}} \leq \mu  \} \subseteq
\{ \wh{\bs{\mu}} < z^* \}$.  It follows that $\Pr \{ \mu \geq
\ovl{\bs{\mu}} \} \leq \Pr \{ \wh{\bs{\mu}} < z^* \} \leq
\f{\de}{2}$. In the case of $\Pr \{ \wh{\bs{\mu}} \leq z^* \} =
\f{\de}{2}$, we have $\f{\de}{2} = \Pr \{ \wh{\bs{\mu}} \leq z^* \}
\leq \mscr{G} (z^*, \mu)$.  Since $\mscr{G} (z, \mu)$ is increasing
with respect to $z \in (0, \mu)$, we have $\{ \mu \geq
\ovl{\bs{\mu}} \} \subseteq \{ \mscr{G} (\wh{\bs{\mu}}, \mu) \leq
\f{\de}{2}, \; \wh{\bs{\mu}} \leq \mu  \} \subseteq \{ \mscr{G}
(\wh{\bs{\mu}}, \mu) \leq \mscr{G} (z^*, \mu), \; \wh{\bs{\mu}} \leq
\mu  \} \subseteq \{ \wh{\bs{\mu}} \leq z^* \}$. It follows that
$\Pr \{ \mu \geq \ovl{\bs{\mu}} \} \leq \Pr \{ \wh{\bs{\mu}} \leq
z^* \} = \f{\de}{2}$.

\epf

\beL $\Pr \{ \mu \leq \udl{\bs{\mu}} \} \leq \f{\de}{2}$. \eeL

\bpf

 For $z > \mu$ and $1 \geq z \geq \f{\ga}{n}$, by Theorem 1 and
 Lemma \ref{lem1}, we have {\small \bee \Pr \{  \wh{\bs{\mu}}
\geq z \} & = & \Pr \li \{ \sum_{i = 1}^{\li \lf \ga \sh z \ri \rf}
X_i \geq \ga \ri \} \leq \exp \li ( \li \lf \ga \sh z \ri \rf
\mscr{M}_{\mrm{B}} \li ( \f{\ga}{ \li \lf \ga \sh z \ri \rf }, \mu \ri ) \ri )\\
& = & \exp \li ( \ga \mscr{M}_{\mrm{I}} \li ( \f{\ga}{ \li \lf \ga
\sh z \ri \rf }, \mu \ri ) \ri ) \leq \exp( \ga \mscr{M}_{\mrm{I}} (
z, \mu) ) \eee} where the last inequality is due to $\f{\ga}{ \li
\lf \ga \sh z \ri \rf } \geq z$ and the fact that
$\mscr{M}_{\mrm{I}} (z, \mu)$ is monotonically decreasing with
respect to $z \in (\mu, 1)$.  For $\mu < z < \f{\ga}{n}$, by Theorem
1 and Lemma \ref{lem1}, we have {\small $\Pr \{  \wh{\bs{\mu}} \geq
z \} = \Pr \li \{  \f{ \sum_{i = 1}^{n} X_i  } { n } \geq z \ri \}
\leq \exp( n \mscr{M}_{\mrm{B}} ( z, \mu) )$}.  Therefore, $\Pr \{
\wh{\bs{\mu}} \geq z \} \leq \mscr{H}( z, \mu)$, where
\[
\mscr{H}( z, \mu) = \bec \exp \li ( \ga \mscr{M}_{\mrm{I}} \li ( z, \mu \ri ) \ri )
&  \tx{for} \; 1 \geq z \geq \f{\ga}{n}, \; z > \mu,\\
\exp( n \mscr{M}_{\mrm{B}} ( z, \mu) )  & \tx{for} \; \mu < z <
\f{\ga}{n}. \eec
\]
Let $z^* \in [0, 1]$ be the smallest number such that $\Pr \{
\wh{\bs{\mu}} > z^* \} \leq \f{\de}{2}$. Then, it must be true that
either $\Pr \{ \wh{\bs{\mu}} \geq z^* \} > \f{\de}{2}$ or $\Pr \{
\wh{\bs{\mu}} \geq z^* \} = \f{\de}{2}$. Observing that $\mscr{H}
(z, \mu)$ is monotonically increasing with respect to $\mu \in (0,
z)$, we have {\small \[ \{ \mu \leq \udl{\bs{\mu}} \} = \{ \mu \leq
\udl{\bs{\mu}} \leq \wh{\bs{\mu}}, \; \mbf{k} > 0 \} \subseteq \li
\{ \mscr{H} (\wh{\bs{\mu}}, \mu) \leq \mscr{H} (\wh{\bs{\mu}},
\udl{\bs{\mu}}) = \f{\de}{2}, \; \mu \leq \udl{\bs{\mu}} \leq
\wh{\bs{\mu}} \ri \} \subseteq \li \{ \mscr{H} (\wh{\bs{\mu}}, \mu)
\leq \f{\de}{2}, \; \wh{\bs{\mu}} \geq \mu \ri \}.
\]}
In the case of $\Pr \{ \wh{\bs{\mu}} \geq z^* \} > \f{\de}{2}$, we
have $\f{\de}{2} < \Pr \{ \wh{\bs{\mu}} \geq z^* \} \leq \mscr{H}
(z^*, \mu)$.  Since $\mscr{H} (z, \mu)$ is decreasing with respect
to $z \in (\mu, 1)$, we have $\{ \mu \leq \udl{\bs{\mu}} \}
\subseteq \{ \mscr{H} (\wh{\bs{\mu}}, \mu) \leq \f{\de}{2}, \;
\wh{\bs{\mu}} \geq \mu  \} \subseteq \{ \mscr{H} (\wh{\bs{\mu}},
\mu) < \mscr{H} (z^*, \mu), \; \wh{\bs{\mu}} \geq \mu  \} \subseteq
\{ \wh{\bs{\mu}} > z^* \}$.  It follows that $\Pr \{ \mu \leq
\udl{\bs{\mu}} \} \leq \Pr \{ \wh{\bs{\mu}} > z^* \} \leq
\f{\de}{2}$. In the case of $\Pr \{ \wh{\bs{\mu}} \geq z^* \} =
\f{\de}{2}$, we have $\f{\de}{2} = \Pr \{ \wh{\bs{\mu}} \geq z^* \}
\leq \mscr{H} (z^*, \mu)$.   Since $\mscr{H} (z, \mu)$ is decreasing
with respect to $z \in (\mu, 1)$, we have $\{ \mu \leq
\udl{\bs{\mu}} \} \subseteq  \{ \mscr{H} (\wh{\bs{\mu}}, \mu) \leq
\f{\de}{2}, \; \wh{\bs{\mu}} \geq \mu  \} \subseteq \{ \mscr{H}
(\wh{\bs{\mu}}, \mu) \leq \mscr{H} (z^*, \mu), \; \wh{\bs{\mu}} \geq
\mu  \} \subseteq \{ \wh{\bs{\mu}} \geq z^* \}$.  It follows that
$\Pr \{ \mu \leq \udl{\bs{\mu}} \} \leq \Pr \{ \wh{\bs{\mu}} \geq
z^* \} = \f{\de}{2}$.

\epf

\end{document}